\documentclass[times,sort&compress,3p]{elsarticle}
\journal{Journal of Multivariate Analysis}

\usepackage{amsmath,amsfonts,amssymb,amsthm,booktabs,color,epsfig,graphicx,hyperref,url}

\usepackage{graphicx}
\usepackage{array}

\usepackage{natbib}
\usepackage{comment}

\def\bpsi{\boldsymbol{\psi}}
\def\bxi{\boldsymbol{\xi}}

\def\bSig{\boldsymbol{\Sigma}}

\usepackage[figuresright]{rotating}

\renewcommand{\hat}{\widehat}

\newtheorem{theorem}{Theorem}

\newtheorem{lemma}{Lemma}

\newtheorem{proposition}{Proposition}

\begin{document}

\begin{frontmatter}

\title{Estimating sparse networks with hubs}

\author[A1]{Annaliza McGillivray}
\author[A2]{Abbas Khalili\corref{mycorrespondingauthor}}
\author[A2]{David A. Stephens}

\address[A1]{Department of Mathematics and Statistics, University of Saskatchewan, Saskatoon, Canada}
\address[A2]{Department of Mathematics and Statistics, McGill University, Montr\'eal, Canada}

\cortext[mycorrespondingauthor]
{Corresponding author. Email address: \url{abbas.khalili@mcgill.ca}}

\begin{abstract}

Graphical modelling techniques based on sparse estimation have been applied to infer 
complex networks in many fields, including biology and medicine, engineering, finance and social sciences. 
One structural feature of some of these networks that poses 
a challenge for statistical inference is the presence of a small number 
of strongly interconnected nodes, which are called hubs. 
For example, in microbiome research hubs or microbial taxa
play a significant role in maintaining stability of the microbial community structure.
In this paper, we investigate the problem of estimating 
sparse networks in which there are a few highly connected hub nodes. 
Methods based on $L_1$-regularization have been widely used for performing 
sparse estimation in the graphical modelling context. However, while 
these methods encourage sparsity, they do not take into account structural 
information of the network. We introduce a new method for estimating networks 
with hubs that exploits the ability of (inverse) covariance estimation methods to include structural information 
about the underlying network. Our method is a weighted lasso approach with novel 
row/column sum weights, which we refer to as the hubs weighted graphical lasso. 
A practical advantage of the new method is that it leads to an optimization problem 
that is solved using the 
efficient graphical lasso algorithm that is already implemented in the \verb+R+ package 
\verb+glasso+. We establish large sample properties of the method when the number of 
parameters diverges with the sample size. We then show via simulations that 
the method outperforms competing methods and illustrate its performance with an application to microbiome data. 


\end{abstract}


\end{frontmatter}

\section{Introduction \label{sec:Introduction}}

Over the past decade, fitting graphical models or networks via estimation of large sparse 
covariance and precision matrices has attracted much attention in modern multivariate analysis.
Applications range from biology and medicine to engineering, economics, finance, and social sciences 
\cite{FanLiaoLiu2016}.
To handle data scarcity in estimating large or high-dimensional sparse networks, 
methods based on $L_1$-regularization 
(\cite{MB2006}, \cite{YuanLin2007}, \cite{FriedmanHastieTibshirani2008}) 
are widely used, the most popular being 
the graphical lasso (\verb+glasso+) of \cite{FriedmanHastieTibshirani2008}.  
The \verb+glasso+ estimates the so-called precision matrix $\Theta = \bSig^{-1}$ via maximizing 
an $L_1$-penalized Gaussian log-likelihood, based on a random sample of $p_n$-dimensional
Gaussian random vectors ${\bf X}_1, \ldots, {\bf X}_n$ with mean ${\bf 0}$ and 
covariance matrix $\bSig$ (see Section \ref{sec:PenalizedLikelihoodFramework}). 
Under the Gaussianity assumption on ${\bf X}_i = (X_{i1}, X_{i2}, \ldots, X_{ip_n})^{\top}$, 
a non-zero element $\theta_{jl}$ of $\Theta$ corresponds to an edge between two 
nodes $X_{ij}$ and $X_{il}$ in a graphical model for the data. The $L_1$-penalty 
 is applied to the off-diagonal elements of the presumably sparse precision matrix $\Theta$.
It is known that the \verb+glasso+ produces a sparse estimate of the precision matrix $\Theta$. 
However, since the $L_1$-penalty increases linearly in $|\theta_{jl}|$, the \verb+glasso+ 
also results in biased estimates of the large $\theta_{jl}$.  
To reduce the estimation bias, \cite{LamFan2009} and \cite{ShenPanZhu2012} 
proposed penalized likelihood approaches based on non-convex penalties such 
as smoothly clipped absolute deviation or 
SCAD \citep{FanLi2001} for  sparse  precision  matrix  estimation and studied  
their  theoretical properties; \cite{FanFengWu2009} introduced the graphical adaptive 
lasso \cite{Zou2006} to attenuate the bias problem in the network estimation.

\begin{figure}
\begin{center}
\includegraphics[scale=.55]{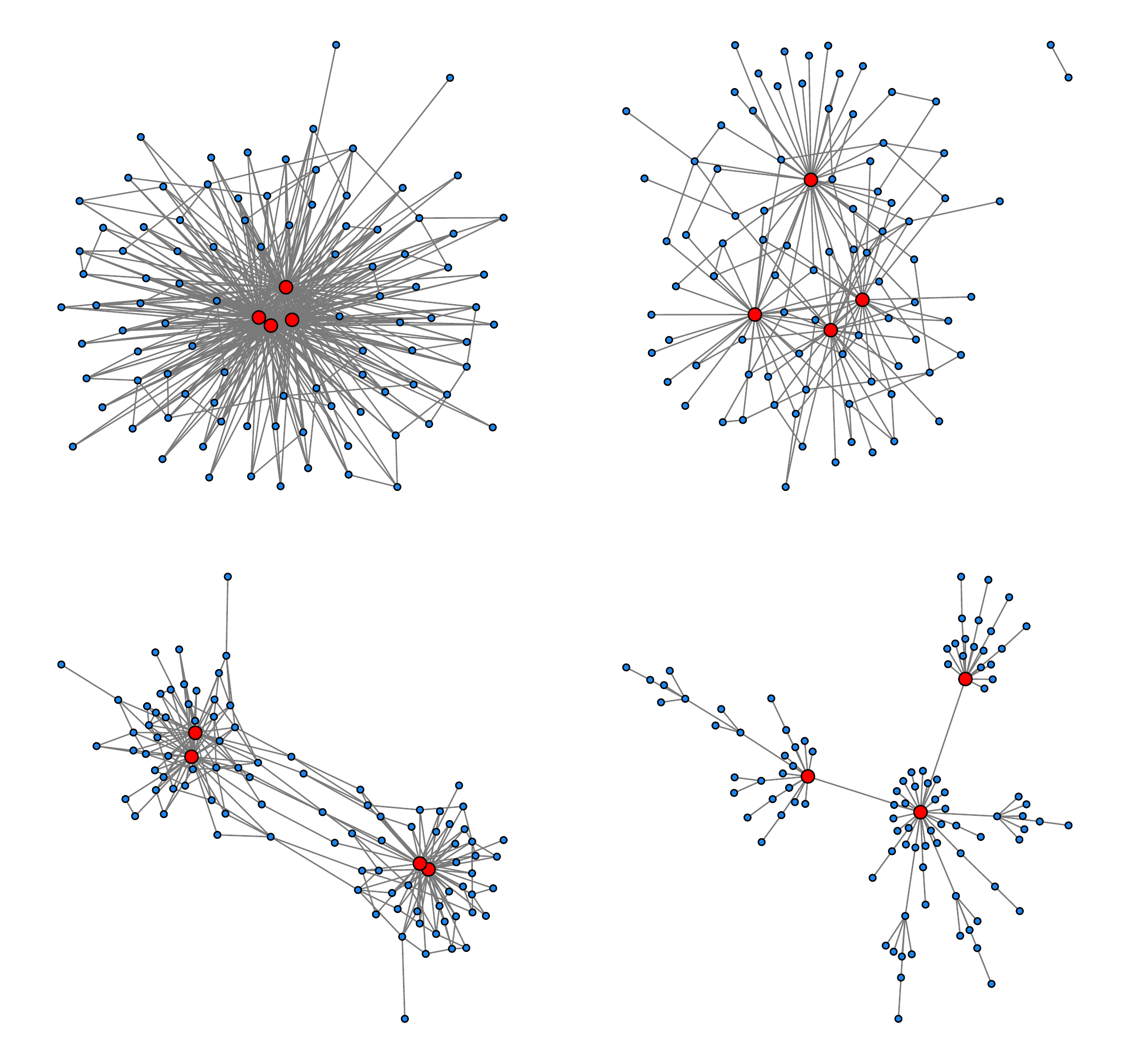}
\end{center}
		\caption{Simulated networks with hub nodes.}
		\label{fig:SimulatedNetworks}
\end{figure}

Penalties such as $L_1$ and SCAD, however, implicitly assume that each potential edge in a network 
is equally likely and/or independent of 
all other edges \citep{Tan2014}, and may thus be inadequate for 
estimating networks with a few highly connected nodes, called {\it stars} or {\it hubs} 
(Figure \ref{fig:SimulatedNetworks}). 
On the other hand, the weights in the graphical adaptive lasso \citep{FanFengWu2009} 
do not take network structural features such as hubs into consideration. 
In this paper, inspired by microbiome data, we propose a 
new regularization method referred to as the hubs weighted graphical lasso 
for estimating sparse networks with a few hub nodes and in the presence of many low-degree nodes. 

Rapidly developing sequencing technologies and analytical techniques have enhanced our ability to study the microorganisms such as bacteria, viruses, archaea and fungi that inhabit the human body \citep{Gilbert2010} 
and a wide range of environments \citep{Turnbaugh2007}. 
The microorganisms inhabiting a particular environment do not exist in isolation, but interact with other microorganisms in a range of mutualistic and antagonistic relationships. 
One goal of microbiome studies is to model these microbial interactions from population-level data as 
a network reflecting co-occurrence and co-exclusion patterns between microbial taxa. This is of interest not only for predicting individual relationships between microbes, but the structure of the interaction networks also gives insight into the organization of complex microbial communities. 
\cite{Friedman2012} used networks of pairwise correlations between microbial taxa to 
model microbe-microbe interactions from microbial abundance data. 
However, correlation can be limiting in the multivariate setting as it is a pairwise measure of 
dependence. In addition, 
statistical challenges in studying networks of microbial interactions arise due to data scarcity 
and the organization of the network's nodes into groups with different levels of connectivity. 
Specifically, microbial association networks tend to be sparse and also 
display hubs \citep{HeiHart2016}.
In ecology, these hubs can represent a few keystone species that are 
vital in maintaining stability of the microbial community \citep{Kurtz2015}.  
 
To accommodate structural information such as hubs 
in network estimation, \cite{Tan2014} proposed the hubs graphical lasso (\verb+HGL+), which is 
a penalization method that encourages 
estimates of the form $\widehat \Theta_n = {\bf Z} + {\bf V} + {\bf V}^\top$, 
where ${\bf Z}$ is a sparse symmetric matrix capturing edges between non-hub 
nodes and ${\bf V}$ is a matrix whose columns are either entirely zero or almost 
entirely non-zero with the non-zero elements of ${\bf V}$ representing hub edges. 
The \verb+HGL+ applies an $L_1$-penalty to the off-diagonal elements of ${\bf Z}$, 
and $L_1$ and group lasso \citep{YuanLin2007} penalties  to the columns of ${\bf V}$. 
The method requires considerable tuning with three tuning parameters present in the $L_1$-penalized likelihood, which are selected by a BIC-type quantity; more details are given in Section \ref{sec:SimulationStudies}. 
The \verb+HGL+ is specifically designed for networks with dense hub nodes, 
referred to as {\it super hubs}. \cite{LiuIhler2011} proposed a method for estimating 
scale-free networks, which are characterized as having a degree distribution 
that follows a power law. 
Their method is, in particular, a re-weighted $L_1$-regularization approach, where the 
weights in the iterative procedure are updated in a rich-get-richer fashion, 
mimicking the generating mechanism of scale-free networks \citep{Barabasi1999}. 
Such an approach, however, cannot model super hubs \citep{Tan2014}. \cite{HeroRaj2012}
 proposed a screening method for hub screening in high dimensions but it does 
 not estimate the edges of the network. 

In this paper, we introduce a new approach for estimating networks with hubs that exploits the ability of 
(inverse) covariance estimation methods to include structural information about the underlying network, and can accommodate both networks with so-called super hubs as well as scale-free networks. More specifically, our method called the hubs weighted graphical lasso (\verb+hw.glasso+), is a weighted graphical lasso approach with novel informative row/column sum weights that allow for differential penalization of hub edges compared to non-hub edges. 
In our theoretical development, we first investigate estimation and selection consistency of a general weighted graphical lasso estimator \citep{FanFengWu2009}, when $p_n \to \infty$ as $n \to \infty$. 
We then provide conditions under which the \verb+hw.glasso+ estimator, based on our proposed weights, achieves the aforementioned theoretical properties. To the best of our knowledge, theoretical properties of \verb+HGL+ \citep{Tan2014} and the method of \citep{LiuIhler2011} are not known. 
In practice, the \verb+hw.glasso+ leads to an optimization problem that 
is solved using the efficient graphical lasso algorithm of \cite{FriedmanHastieTibshirani2008}, 
already implemented in the \verb+R+ package \verb+glasso+. 
Via extensive simulations we show that, in comparison to competing methods, 
the \verb+hw.glasso+ performs well in finite sample situations considered here. 

The remainder of this paper is organized as follows. In Section \ref{sec:PenalizedLikelihoodFramework}, we introduce the penalized likelihood problem and commonly used penalty functions in the context of performing sparse inverse covariance estimation. In Section \ref{sec:ProposalNetworkSelection}, we present the hubs weighted graphical lasso (\verb+hw.glasso+) estimator, and investigate its theoretical properties in Section \ref{sec:AsymptoticProperties}. We then assess its finite sample performance through simulation studies in Section \ref{sec:SimulationStudies}, and with an application to two microbiome data sets in Section \ref{sec:MicrobiomeApplication}. We conclude with a discussion in Section \ref{sec:Conlusions}. 
The proofs of our main theoretical results are provided in the Appendix. Additional simulations are given in our Supplementary Material. \\

\section{Problem Setup \label{sec:PenalizedLikelihoodFramework}}

Suppose ${\bf X}_1, \ldots, {\bf X}_n$ are $p_n$-dimensional independent and identically distributed (iid) random vectors from a Gaussian distribution with mean ${\bf 0}$ and covariance matrix $\bSig = \Theta^{-1}$, and let ${\bf x}_1, \ldots, {\bf x}_n$, with ${\bf x}_i = (x_{i1}, \ldots, x_{ip_n})^\top$, denote realizations of the random variables. Further, denote the sample covariance matrix by 
$S_n$, where $S_n=\sum_{i=1}^n {\bf x}_i {\bf x}_i^{\top}/n$. Then the re-scaled log-likelihood function of $\Theta$ (up to a constant) is given by
\begin{align}
\label{eq:like}
\ell_n(\Theta) = \log{\det{(\Theta)}} - \text{tr}(S_n\Theta),
\end{align}
where $\det(\cdot)$ and $\text{tr}(\cdot)$ respectively denote the determinant and trace. 
For a sparse Gaussian graphical model, the precision matrix $\Theta$ is estimated by the maximizer of the penalized log-likelihood function
\begin{equation}
\label{eq:penlike}
pl_n(\Theta) = \ell_n(\Theta) - p_{\lambda_n}(\Theta), 
\end{equation}
where $p_{\lambda_n}(\cdot)$ is a generic penalty function on $\Theta$ with tuning parameter $\lambda_n>0$. 

\cite{FriedmanHastieTibshirani2008} considered the $L_1$-penalty function $p_{\lambda_n}(\Theta)=\lambda_n \sum_{i < j} |\theta_{ij}|$ in \eqref{eq:penlike} and proposed the graphical lasso (\verb+glasso+) algorithm that makes use of a block coordinate descent procedure to optimize \eqref{eq:penlike}. \cite{LamFan2009} studied nonconvex penalties such as the SCAD in \eqref{eq:penlike}. 
While these penalties induce sparsity in the estimated $\Theta$, 
they do so by penalizing the elements of $\Theta$ equally and/or independently of each other. One penalty function that allows for varying levels of penalization to the entries $\theta_{ij}$ is the adaptive lasso (\citep{FanFengWu2009}, \citep{Zou2006}), given by $p_{\lambda_n}(\Theta)= \lambda_n \sum_{i < j} \tilde{w}_{ij}|\theta_{ij}|$, where $\tilde{w}_{ij} = 1/ |\tilde{\theta}_{ij}|^{\gamma}$ for some $\gamma>0$ and any consistent estimate $\widetilde{\Theta}_n=(\tilde{\theta}_{ij})_{i,j=1}^{p_n}$ of $\Theta$. While these choices of the penalty result in a sparse estimate of $\Theta$ and lead to desirable asymptotic properties \citep{LamFan2009}, they do not incorporate any prior information of structural features such as hub nodes in the precision matrix. This motivated us to propose a method that allows for the inclusion of such rich structural information
in the penalty function in \eqref{eq:penlike}, and our numerical study shows that 
this consideration greatly enhances finite sample performance of the method. 

\section{Hubs Weighted Graphical Lasso \label{sec:ProposalNetworkSelection}}
 
 When the true underlying Gaussian graphical model has hub nodes (see Figure \ref{fig:SimulatedNetworks}), 
we wish to incorporate into our estimation procedure not merely sparsity but the knowledge of 
the presence of these highly connected nodes. In this section, we present a new penalty 
function in (\ref{eq:penlike}) that utilizes this knowledge.

 Since the true underlying graph has hub nodes, in the precision matrix $\Theta$
 the rows/columns corresponding to each hub node are significantly denser (i.e., have more non-zero elements) than those corresponding to the non-hub nodes. In Figure \ref{fig:SimulatedNetworks}, we display four different types of networks with hubs, from top-left to bottom-right: the first, illustrates a network with so-called ``super hubs'', while 
 the second and third display networks with hubs that are less densely connected than the 
 ``super hubs'' in the first, and the fourth displays a scale-free network \citep{Barabasi1999}. 

To estimate networks with hubs, we adopt a new weighted lasso approach that uses informative weights based on row/column sums of the precision matrix. In what follows, we outline our proposed estimation procedure by first introducing the new weights. 

 Let $\widetilde{\Theta}_n=(\tilde{\theta}_{ij})_{i,j=1}^{p_n}$ be any consistent estimator of the precision matrix 
 $\Theta_0$. 
 We may take $\widetilde{\Theta}_n$ to be the 
 precision matrix estimator obtained from the graphical lasso
 \citep{FriedmanHastieTibshirani2008}, which is consistent under the conditions of Theorem 1 below \citep{Rothman2008}.
 We then construct the symmetric matrix $\widetilde{W}_n=(\widetilde{w}_{ij})_{i,j=1}^{p_n}$ of weights 
\begin{align}
\widetilde{w}_{ij} =  
      \displaystyle \frac{1}{  \displaystyle | \tilde{\theta}_{ij} |^{\gamma_1} \left(  \| \tilde{\boldsymbol{\theta}}_{\neg i}  \|_1 \cdot   \| \tilde{\boldsymbol{\theta}}_{\neg j} \|_1 \right)^{\gamma_2}} \hspace{0.5cm}, \text{ if } i \neq j       
       \label{eq:Weights}
\end{align}
and $\widetilde{w}_{ij}=0$ if $i=j$, for some $\gamma_1$, $\gamma_2>0$, where $\tilde{\boldsymbol{\theta}}_{\neg i}= \left\{ \tilde{\theta}_{ik}:  k=1, \ldots, p_n, k \neq i  \right\}$ is the $i^{\text{th}}$ row (or by symmetry, the $i^{\text{th}}$ column) of $\widetilde{\Theta}_n$, and $\| \tilde{\boldsymbol{\theta}}_{\neg i} \|_1=\sum_{k \neq i}  |\tilde{\theta}_{ik}|$.

We now define the hubs weighted graphical lasso (\verb+hw.glasso+) estimator $\widehat{\Theta}_n$ of $\Theta$ to be 
\begin{align}
\widehat{\Theta}_n = \displaystyle \operatorname*{arg\,max}_{\Theta \succ 0} \left\{  \ell_n(\Theta) - \lambda_n \| \widetilde{W}_n * \Theta \|_{1}  \right\},
\label{eq:HWGL1}
\end{align}
where $\lambda_n>0$ is a tuning parameter, and $*$ is the Schur matrix product so that 
\begin{align}
\label{eq:pen1-function}
\| \widetilde{W}_n* \Theta \|_{1}  = \sum_{i<j} \widetilde{w}_{ij}|\theta_{ij}|.
\end{align}

 The proposed approach belongs to the family of weighted graphical lasso methods that allow for 
different penalties on the entries of $\Theta$, which includes the graphical adaptive lasso 
\citep{FanFengWu2009}. Weighted lasso approaches can result in less bias 
than the standard lasso by adapting penalties to incorporate information 
about the location of zeros, based on either an initial estimate or background knowledge.
 
 The weights $\widetilde{w}_{ij}$ in (\ref{eq:Weights}) are designed to allow for less penalization of hub 
 edges compared to non-hub edges. If $\theta_{ij}^0=0$, similar to the adaptive lasso \citep{FanFengWu2009}, the weights are expected to get inflated (to infinity as the sample size grows) because of the presence of the term $| \tilde{\theta}_{ij} |^{\gamma_1}$, irrespective of whether $i$ and $j$ are hubs. For $\theta_{ij}^0 \neq 0$ with at least one of $i$ and $j$ hubs, given the signal strength assumption (Condition 2 in Section \ref{sec:AsymptoticProperties}), we expect $(\| \tilde{\boldsymbol{\theta}}_{\neg i}  \|_1 \cdot   \| \tilde{\boldsymbol{\theta}}_{\neg j} \|_1)$ to be large (greater than 1) due to the hub structure, which results in smaller weights $\widetilde{w}_{ij}$ in (\ref{eq:Weights}) compared to the weights $|\tilde{\theta}_{ij}|^{-\gamma_1}$ in the standard adaptive lasso \citep{FanFengWu2009}. For $\theta_{ij}^0 \neq 0$ with neither $i$ nor $j$ hubs, then it is expected that the proposed \verb+hw.glasso+ performs similarly to the adaptive lasso. 

It is worth mentioning that the role of the penalty function (\ref{eq:pen1-function}) is not to do a group selection, where the group would correspond to the hub rows/columns, but rather to allow for different levels of penalization on $\theta_{ij}$ based on an initial consistent estimator $\widetilde{\Theta}_n$. This is in contrast to the penalty function in \citep{Mohan2014}, which is a group lasso penalty applied to the rows/columns. In this case, an overlap issue arises since the $(i,j)^\text{th}$ entry of the matrix is contained in both the $i^\text{th}$ and $j^\text{th}$ groups. As a result, the group lasso penalty with overlapping groups no longer selects groups (i.e., leaving hub rows/columns fully non-zero). The group lasso penalty with overlapping groups in the regression context is also discussed in \citep{Obozinski2011}.\\
	
\noindent {\bf Numerical Algorithm:}\\

The advantage of the \verb+hw.glasso+ method is that it leads to the optimization problem \eqref{eq:HWGL1} that can be solved using the efficient graphical lasso algorithm of \cite{FriedmanHastieTibshirani2008}, already implemented in the \verb+R+ package 
\verb+glasso+. In their implementation, the user may specify a symmetric weight matrix, 
which in our case is $\widetilde{W}_n$ defined in \eqref{eq:Weights}. 
For the choice of the tuning parameter $\lambda_n$, we employ the Bayesian information criterion 
(BIC) which has been widely used in the literature \cite{Gao2012}. In our simulation 
studies and real data analysis, respectively, in Sections \ref{sec:SimulationStudies} and \ref{sec:MicrobiomeApplication}, we take $\gamma_1=\gamma_2=1$. 
 
\section{Theoretical Properties \label{sec:AsymptoticProperties}}
 
In this section, we first view the estimator $\widehat{\Theta}_n$ in (\ref{eq:HWGL1}) as a general 
weighted \verb+glasso+ estimator and derive conditions on the weights $\widetilde{w}_{ij}$ 
in \eqref{eq:pen1-function} that  guarantee consistency and sparsistency (see below) 
of $\widehat{\Theta}_n$. We then focus on the specific weights \eqref{eq:Weights} 
that resulted in our hubs estimator 
\verb+hw.glasso+. 
The weights $\widetilde{w}_{ij}$ typically depend
on the sample size $n$ and are possibly random. 

We assume that ${\bf X}_1, \ldots, {\bf X}_n$ are $p_n$-dimensional iid Gaussian random vectors with mean ${\bf 0}$ and true covariance matrix $\bSig_0$. 
The corresponding true sparse precision matrix is $\bSig_0^{-1}=\Theta_0 = (\theta_{ij}^{0})_{i,j=1}^{p_n}$, where  
$p_n \to \infty$ at a certain rate to be later specified, as $n \rightarrow \infty$. 
First, we introduce some notation and state certain regularity conditions on the true precision matrix 
$\Theta_0$.

We define $T=\left\{ (i,j): \theta_{ij}^{0} \neq 0, i<j \right\} \neq \emptyset$ to be the set of indices of all non-zero off-diagonal elements in $\Theta_0$ and let $q_n = | T |$ be the cardinality of $T$. The set of indices of the true zero elements of $\Theta_0$ is denoted by $T^c$. Let $\phi_{\text{min}}(A)$ and $\phi_{\text{max}}(A)$ denote the minimum and maximum eigenvalues of a matrix $A$. Further, 
let $\| A \|_F^2 = \text{tr}(A^\top A)$ and 
$\| A \|^2 = \phi_{\text{max}}(A^\top A)$ be the Frobenius and operator norms of 
$A$, respectively. Also, recall from \eqref{eq:Weights} that $\|\boldsymbol{\theta}^{0}_{\neg i}  \|_1=\sum_{k \neq i} |\theta_{ik}^{0}|$. We assume that the following regularity conditions hold. \\

\noindent \emph{Condition 1:} There exist constants $\tau_1$ and $\tau_2$ such that $0 < \tau_1 \le \phi_{\text{min}}(\Theta_0) < \phi_{\text{max}}(\Theta_0) \le \tau_2 < \infty$.

\noindent \emph{Condition 2:} There exists a constant $\tau_3>0$ such that $\min_{(i,j) \in T} | \theta_{ij}^{0} | \geq \tau_3$. \\

Condition 1 guarantees the existence of the true inverse covariance matrix $\Theta_0$, which is required under the Gaussianity assumption. If this condition is violated, the Gaussian model 
may no longer be appropriate for this problem, which then calls for alternative models. 
Condition 2 is a signal strength assumption; it ensures that the non-zero elements of $\Theta_0$ are bounded 
away from zero. The proofs of our results are given in the Appendix. 
Our first result concerns the estimation {\it consistency} of the weighted 
\verb+glasso+ estimator. 

\begin{theorem}
\label{thm:consistency}
(Consistency) 
Suppose Conditions 1 and 2 hold, and $(p_n+q_n)(\log{p_n})/n=o(1)$. 
Further, assume that $\lambda_n$ and $\widetilde{w}_{ij}$ are chosen such that 
$\lambda_n \max_{(i,j) \in T} \widetilde{w}_{ij} = O_p( \left\{ (\log{p_n})/n \right\}^{1/2} )$ and 
$\left\{ (\log{p_n})/n \right\}^{1/2} \left\{ \min_{(i,j) \in T^c} \widetilde{w}_{ij}\right\}^{-1} = O_p(\lambda_n)$. 
Then the weighted \verb+glasso+ estimator $\widehat{\Theta}_n$ satisfies 
\begin{align}
\| \widehat{\Theta}_n - \Theta_0 \|_{\text F} = O_p\left[ \left\{  \frac{(p_n+q_n)\log{p_n}}{n} \right\}^{1/2} \right].
\label{eq:RateOfConvergence}
\end{align}
\end{theorem}

Theorem \ref{thm:consistency} shows that with the proper choice of the tuning 
parameter $\lambda_n$ and the weights $\widetilde{w}_{ij}$, 
$\widehat{\Theta}_n$ is a consistent estimator of $\Theta_0$. For example, 
if the (possibly random) weights are chosen such that  
$\max_{(i,j) \in T} \widetilde{w}_{ij} \rightarrow C_1<\infty$ 
and $\min_{(i,j) \in T^c} \widetilde{w}_{ij} \rightarrow \infty$ (in probability), as $n \to \infty$, then 
the choice $\lambda_n = C_2 \left\{  (\log{p_n})/n \right\}^{1/2}$, 
for some finite constant $C_2>0$, results in consistency of  
$\widehat{\Theta}_n$.
As pointed out by \cite{LamFan2009} and \cite{Rothman2008}, 
the worst part of the rate of convergence of $\widehat{\Theta}_n$ in 
(\ref{eq:RateOfConvergence}) is the term $p_n \log{p_n}/n$, which is 
due to the estimation of $p_n$ diagonal elements of $\Theta_0$. This rate can be improved to 
$\left\{(q_n \log{p_n})/n \right\}^{1/2}$ if we were to estimate the inverse of the true correlation matrix; 
more details are given in Remark 1 below. The effect of diverging 
dimensionality is reflected by the term $\log{p_n}$. 

Our next result establishes that the weighted 
\verb+glasso+ 
estimates the true zero entries of the precision matrix as zero 
with probability tending to 1. This property is referred to as 
\emph{sparsistency} in \cite{LamFan2009}.

\begin{theorem}
\label{thm:sparsistency}
(Sparsistency) Assume the conditions of Theorem 1 are fulfilled, and that
\begin{align}
\label{sparse-rate}
\| \widehat{\Theta}_n - \Theta_0 \|^2 &= O_p(\eta_n)
\end{align}
for a sequence $\eta_n$ such that $\eta_n \rightarrow 0$ and 
$\left\{ \min_{(i,j) \in T^c} \widetilde{w}_{ij}\right\}^{-2} \eta_n  = O_p(\lambda_n^2)$. 
Then the weighted \verb+glasso+ estimator $\widehat{\Theta}_n$ 
has the property $P(\hat{\theta}_{ij}=0: (i,j) \in T^c) \rightarrow 1$, as $n \rightarrow \infty$.
\end{theorem}

The two theorems provide general conditions on the weights $\widetilde{w}_{ij}$ and
$(\lambda_n, \eta_n)$ that guarantee consistency and sparsistency  
of the weighted \verb+glasso+ estimator. In this paper, we focus on the specific 
weights $\widetilde{w}_{ij}$ in \eqref{eq:Weights} which are used in our hubs weighted graphical 
lasso (\verb+hw.glasso+) estimator. These weights are constructed 
using the popular \verb+glasso+ estimator $\widetilde \Theta_n$, which is a consistent estimator of $\Theta_0$ \citep{Rothman2008} under the conditions of Theorem 1 above.
Proposition \ref{prop1} below verifies conditions of the theorems on such weights. 

\begin{proposition}
\label{prop1}{\it 
Consider the \verb+hw.glasso+ estimator $\widehat \Theta_n$ with the specific 
weights $\widetilde w_{ij}$ in \eqref{eq:Weights}. 

\begin{itemize}
\item[{\bf (a)}] If there exists a pair $(i,j) \in T^c$ such that $\|\boldsymbol{\theta}^{0}_{\neg i}  \|_1 \neq 0$ and $\|\boldsymbol{\theta}^{0}_{\neg j}  \|_1 \neq 0$, then the estimator has consistency 
property \eqref{eq:RateOfConvergence} if  
$\{\lambda_n, p_n, q_n\}$ satisfy 
\begin{align}
\label{final-condition1a}
\lambda_n = O(\sqrt{\log p_n/n})~,~
\left( \frac{\log{p_n}}{n} \right)^{1/2} 
\left\{ \frac{(p_n+q_n)\log{p_n}}{n} \right\}^{\gamma_1/2} \lambda^{-1}_n 
& = O(1).
\end{align}
The estimator has also sparsistency property if we have \eqref{final-condition1a} and 
\eqref{sparse-rate} with 
$\eta_n$ satisfying 
\begin{align}
\label{final-condition2a}
\sqrt \eta_n \left( \frac{(p_n+q_n)\log{p_n}}{n} \right)^{\gamma_1/2} \lambda_n^{-1} 
& = O(1).
\end{align} 
\item[{\bf (b)}] If for all $(i,j) \in T^c$, $\|\boldsymbol{\theta}^{0}_{\neg i}  \|_1 = 0$ or $\|\boldsymbol{\theta}^{0}_{\neg j}  \|_1 = 0$, then the estimator has consistency 
property \eqref{eq:RateOfConvergence} if  
$\{\lambda_n, p_n, q_n\}$ satisfy 
\begin{align}
\label{final-condition1b}
\lambda_n = O(\sqrt{\log p_n/n})~,~
\left( \frac{\log{p_n}}{n} \right)^{1/2} 
\left\{ \frac{(p_n+q_n)\log{p_n}}{n} \right\}^{\gamma_1/2} 
\left\{ \frac{p_n (p_n+q_n)\log{p_n}}{n}
\right\}^{\gamma_2/2}
\lambda^{-1}_n 
& = O(1).
\end{align}
The estimator has also sparsistency property if we have \eqref{final-condition1b} and 
\eqref{sparse-rate} with 
$\eta_n$ satisfying 
\begin{align}
\label{final-condition2b}
\sqrt \eta_n \left( \frac{(p_n+q_n)\log{p_n}}{n} \right)^{\gamma_1/2} 
\left\{ \frac{p_n (p_n+q_n)\log{p_n}}{n} \right\}^{\gamma_2/2}
\lambda_n^{-1} 
& = O(1).
\end{align} 
\end{itemize}
}
\end{proposition}

As per the conditions of Theorems \ref{thm:consistency} and \ref{thm:sparsistency}, 
the quantity 
$
\left\{ \min_{(i,j) \in T^c} \widetilde{w}_{ij}\right\}^{-1} 
= 
\max_{(i,j) \in T^c} \left\{  | \tilde{\theta}_{ij}|^{\gamma_1} 
  \Big[ \| \tilde{\boldsymbol{\theta}}_{\neg i}  \|_1 \| 
  \tilde{\boldsymbol{\theta}}_{\neg j}  \|_1 \Big]^{\gamma_2} \right\}
$ 
plays an important role in the properties of the proposed estimator $\widetilde \Theta_n$.
Under case (a) of the proposition, we have that 
$\max_{(i,j) \in T^c} (\|\boldsymbol{\theta}^{0}_{\neg i} \|_1 
\|\boldsymbol{\theta}^{0}_{\neg j}\|_1) \not= 0$. Due to the consistency of 
the initial estimator $\widetilde \Theta_n$ of $\Theta_0$, the quantity 
$\max_{(i,j) \in T^c} \left\{ 
  \Big[ \| \tilde{\boldsymbol{\theta}}_{\neg i}  \|_1 \| 
  \tilde{\boldsymbol{\theta}}_{\neg j}  \|_1 \Big]^{\gamma_2} \right\}
$
converges to a non-zero value, in probability, as $n \to \infty$. 
Thus, it is not surprising that in this case, our proposed weights in \eqref{eq:Weights} 
asymptotically behave similar to the standard weights $|\tilde \theta_{ij}|^{-\gamma_1}$ 
in the graphical adaptive lasso \citep{FanFengWu2009} estimator, i.e. 
the weights in \eqref{eq:Weights} with $\gamma_1>0$ and $\gamma_2 = 0$.
On the other hand, under case (b) of the proposition, we have that 
$\max_{(i,j) \in T^c} (\|\boldsymbol{\theta}^{0}_{\neg i} \|_1 
\|\boldsymbol{\theta}^{0}_{\neg j}\|_1) = 0$. Thus, due to the consistency of 
the initial estimator $\widetilde \Theta_n$ of $\Theta_0$, 
the quantity 
$\max_{(i,j) \in T^c} \left\{ 
  \Big[ \| \tilde{\boldsymbol{\theta}}_{\neg i}  \|_1 \| 
  \tilde{\boldsymbol{\theta}}_{\neg j}  \|_1 \Big]^{\gamma_2} \right\}
$
converges to zero, in probability, as $n \to \infty$. Therefore, in this case both tuning parameters 
$(\gamma_1, \gamma_2)$ play a role in the behaviour of our proposed estimator. 
In our simulation study in Section \ref{sec:SimulationStudies}, 
we have also examined the effects of these tuning parameters on the
finite sample performance of our proposed \verb+hw.glasso+ estimator compared to its competitors.


We now discuss the rate $\eta_n$ in \eqref{sparse-rate}. 
Note that for any $m \times m$ matrix $A$, 
we have that $\|A \|^2 \le \| A \|^2_{\text F} \le m~ \| A \|^2$. 
Let $r_n = \sqrt{(p_n+q_n)\log{p_n}/n}$. Under \eqref{final-condition1a} or \eqref{final-condition1b}, 
the \verb+hw.glasso+ estimator $\widehat \Theta_n$ satisfies 
\eqref{eq:RateOfConvergence} and thus 
$\|\widehat \Theta_n  - \Theta_0\|^2 \le r^2_n \le p_n \|\widehat \Theta_n  - \Theta_0\|^2 $. 
If we consider the worst case scenario that $\eta_n = r^2_n$, 
then the sparsistency conditions \eqref{final-condition2a} and 
\eqref{final-condition2b}, respectively, become
\begin{equation}
\label{dim-con1}
\left\{ (p_n+q_n)^{\gamma_1+1} \left ( \frac{\log{p_n}}{n} \right )^{\gamma_1} \right\}^{1/2} = O(1)
~~\text{and}~~~
\left\{ (p_n+q_n)^{\gamma_1+1} \left ( \frac{\log{p_n}}{n} \right )^{\gamma_1} \right\}^{1/2}
\left\{ \frac{p_n (p_n+q_n)\log{p_n}}{n} \right\}^{\gamma_2/2}
= O(1).
\end{equation}

On the other hand, in the optimistic scenario that  
$\eta_n = r^2_n/p_n$, conditions 
\eqref{final-condition2a} and \eqref{final-condition2b}, respectively, become
\begin{equation}
\label{dim-con2}
\sqrt{\frac{(p_n+q_n)}{p_n}} \left( \frac{(p_n+q_n)\log{p_n}}{n} \right)^{\gamma_1/2} = O(1)
~~\text{and}~~~
\sqrt{\frac{(p_n+q_n)}{p_n}} \left( \frac{(p_n+q_n)\log{p_n}}{n} \right)^{\gamma_1/2} 
\left\{ \frac{p_n (p_n+q_n)\log{p_n}}{n} \right\}^{\gamma_2/2}
= O(1).
\end{equation}

Thus, under the above two scenarios considered for $\eta_n$, 
as long as \eqref{dim-con1} or \eqref{dim-con2} are satisfied, 
the \verb+hw.glasso+ estimator $\widehat \Theta_n$ has the sparsistency property. \\

\noindent
 {\bf Remark 1:} As mentioned above, the worst part of the rate of convergence of $\widehat{\Theta}_n$ in (\ref{eq:RateOfConvergence}) is $p_n \log{p_n}/n$ because of the estimation of $p_n$ diagonal elements of $\Theta_0$. It turns out that the rate can be improved as follows. Using the sample correlation matrix $R_n$ in (\ref{eq:like}) instead of the sample covariance matrix $S_n$, we obtain the penalized estimator, say $\widehat{K}_n$ of the true inverse correlation matrix $R_0^{-1}$, solving a similar optimization problem as in (\ref{eq:HWGL1}). Here the weights are constructed based on a consistent estimator of the inverse correlation matrix. We then define a modified correlation-based estimator of $\Theta_0$ by $\widehat{\Theta}_n^* = \widehat{D}^{-1} \widehat{K}_n \widehat{D}^{-1}$, where $\widehat{D}$ is the diagonal matrix of the sample standard deviations. Similar to Theorem 2 of \citep{Rothman2008} and Theorem 3 of \citep{LamFan2009}, we obtain the rate of convergence of $\widehat{\Theta}_n^*$ to $\Theta_0$ in terms of the operator norm,
 $\| \widehat{\Theta}_n^* - \Theta_0 \| = O_P(\eta_n)$, where $\eta_n = \left\{ (1+q_n) \log{p_n}/n  \right\}^{1/2}$.

\section{Numerical Results \label{sec:SimulationStudies}}

In this section, we compare via simulation the finite sample performance of our proposed 
\verb+hw.glasso+ procedure to the graphical lasso 
(\verb+glasso+, \cite{FriedmanHastieTibshirani2008}), 
the graphical adaptive lasso (\verb+Ada-glasso+) \citep{FanFengWu2009}, the scale-free (\verb+SF+) 
network estimation procedure of \cite{LiuIhler2011}, and the hubs graphical lasso (\verb+HGL+) of 
\cite{Tan2014}. We also provide simulation results for a two-step \verb+hw.glasso+ procedure, 
introduced in Section \ref{sec:TwoStepApproach}, in the case where the hubs are \emph{unknown}, 
but also in the case where the hubs are \emph{known} which is a reasonable assumption in 
some biological applications. 

In practice, both \verb+hw.glasso+ and \verb+Ada-glasso+ require an initial estimator to construct their corresponding weights. In our simulations, we considered three choices of such an initial estimator: the inverse of the sample covariance matrix $S_n$ when $n>p$, and the \verb+glasso+ estimator and the inverse of the shrunken sample covariance matrix $S_n + \alpha I_p$, for some $\alpha>0$, in both cases $n>p$ and $n \leq p$. When $n > p$, all three choices yield similar results, but when $n \leq p$, the inverse of the shrunken sample covariance matrix yielded better results, which are reported in Tables \ref{tab:Simulation1} to \ref{tab:Simulation4}. 

 To implement the graphical (adaptive) lasso and our method, we use the \verb+R+ function \verb+glasso+ and 
 select the tuning parameter $\lambda$ from a fine grid based on BIC, and we 
 set $\gamma_1 = \gamma_2 = 1$. To implement \verb+HGL+, 
 we use the \verb+R+ package \verb+hglasso+. The \verb+HGL+ requires the selection of three tuning parameters $\rho_1, \rho_2$ and $\rho_3$, along with a user-specified parameter $c$ in a BIC-type quantity in \citep{Tan2014} that is used to select the $\rho_i$'s from fine grids. We consider different values of $c$ in our simulations.

 \subsection{Performance Measures and Simulation Settings \label{sec:PerformanceMeasuresSettings}}

 We now provide the performance measures by which various procedures are assessed as well as the simulation settings under consideration. We first introduce some notation. Let TP, TN, FP and FN denote the numbers of true positives (true non-zero $\theta^0_{ij}$'s), true negatives (true zero $\theta^0_{ij}$'s), false positives, and false negatives, respectively. Further, let $\mathcal{H}$ denote the set of indices of true hub nodes, $\widehat{\mathcal{H}}$ the set of indices of  estimated hub nodes, and $| \mathcal{H} |$ denote the size of the set $\mathcal{H}$. To assess the hub structure recovery performance of each of the methods, we consider a node to be a hub if it is connected to more than $k\%$, for some $k$, 
 of all other nodes. 
The methods are evaluated using the following empirical measures: 
\begin{itemize}
	\item True negative rate (TNR, specificity):	 $\displaystyle
		\text{TNR}= \displaystyle \frac{\text{TN}}{\text{TN} + \text{FP}} = \frac{\sum_{\substack{i < j}} I\left( \hat \theta_{ij}=0, \theta_{ij}^{0}=0 \right)} {\sum_{i < j} I\left( \theta_{ij}^{0}=0 \right)} $
\item True positive rate (TPR, sensitivity):
		$\text{TPR}=\displaystyle \frac{\text{TP}}{\text{TP} + \text{FN}} =  \frac{\sum_{i \leq j} I\left( \hat{\theta}_{ij} \neq 0, \theta_{ij}^{0} \neq 0 \right)}{\sum_{i \leq j} I\left( \theta_{ij}^{0} \neq 0 \right)}$
	\item Percentage of correctly estimated hub edges: $\displaystyle\frac{\sum_{i \in \mathcal{H}, i \neq j} I\left( \hat{\theta}_{ij} \neq 0, \theta_{ij}^{0} \neq 0 \right)}{ \sum_{i \in \mathcal{H}, i \neq j} I \left( \theta_{ij}^{0} \neq 0 \right)} \times 100\% $

	\item Percentage of correctly estimated hub nodes: $\displaystyle \frac{|\widehat{\mathcal{H}} \cap \mathcal{H} |}{|\mathcal{H}|} \times 100\%$
	\item Percentage of correctly estimated non-hub nodes:
		$\displaystyle \frac{|\widehat{\mathcal{H}}^c \cap \mathcal{H}^c |}{|\mathcal{H}^c|} \times 100\%$, where $\widehat{\mathcal{H}}^c = \left\{1, \ldots, p \right\} \backslash {\mathcal{H}}$

	\item Frobenius norm: $\frac{1}{p} \| \widehat{\Theta}_n - \Theta_{0} \|^2_{\text F} = \frac{1}{p} \sum_{i \neq j} (\hat{\theta}_{ij} - \theta_{ij}^{0})^2$, 
	
\end{itemize}
where $\widehat{\Theta}_n =(\hat{\theta}_{ij})_{i,j=1}^p$ is the estimated precision matrix, and $\Theta_{0}=(\theta_{ij}^{0})_{i,j=1}^p$ represents the true underlying precision matrix (network). Averages (and standard errors) of these performance measures over 100 replications are reported in Tables \ref{tab:Simulation1} to \ref{tab:Simulation4}. 

We consider four generating mechanisms for the adjacency matrix $A$ of the network, similar to those in \cite{Tan2014}: 
\begin{enumerate}
	\item[(i)]  We randomly select the set $\mathcal{H}$ of true hub nodes and set the elements of the corresponding rows/columns of the adjacency matrix $A$ equal to 1 with probability 0.8 and 0 otherwise. Next, we set $A_{ij} = A_{ji} = 1$ for all $i <j$ with probability 0.01, and 0 otherwise. 
	
	\item[(ii)] 
	We use the same setup as in (i) except that, to generate the adjacency matrix $A$, each hub node is 
	connected to another node with probability 0.3. 
	
	\item[(iii)] The adjacency matrix is
	\(
		A = \begin{pmatrix} 
			A_{1} & B \\
			B^T & A_2 \\
		\end{pmatrix}, \nonumber
	\)
	where $A_1$ and $A_2$ are generated as in (i), except that all nodes have a connection probability of 0.04, and $B=(b_{ij})$ has $b_{ij}=1$ with probability 0.01 and $b_{ij}=0$ otherwise.
	
\item[(iv)] Scale-free networks: for a scale-free network, the probability that a node has degree $d$ follows a power law distribution $P(d) \sim d^{-\alpha}$. Such a network is generated using the algorithm in \cite{Barabasi1999} that incorporates growth and preferential attachment, which are two mechanisms that are common to a number of real-world networks, such as business networks and social networks. We use the \verb+R+ package \verb+igraph+ to generate scale-free networks with $\alpha=2$. Note that the hub nodes in this simulation are less densely connected than those in Simulations (i)-(iii).
\end{enumerate}

For each of the adjacency matrices in (i)-(iv), we then construct a symmetric matrix $\Omega$ such that $\Omega_{ij}= 0$ if $A_{ij}=0$, and $\Omega_{ij}$ are independent from the uniform distribution on $[-0.8,-0.5] \cup [0.5,0.8]$ if $A_{ij}=1$. Finally, we take $\Theta = \Omega + \left\{ 0.1 - \lambda_{\text{min}}(\Omega)\right\}I_p$, where $\lambda_{\text{min}}(\Omega)$ is the smallest eigenvalue of $\Omega$, to ensure that all eigenvalues of $\Theta$ are positive. For Simulations (i) and (ii), we take the number of hubs to be $| \mathcal{H} |=\lfloor p/25 \rfloor$.  The simulated networks for $p=100$ are displayed in Figure \ref{fig:SimulatedNetworks}. When evaluating the performance of each of the methods, we consider a node to be a hub if it is connected to more than $(k=10)\%$ of all other nodes. Note that for Simulations (i)-(iii), there is a clear distinction between hubs and non-hubs, but the cutoff threshold of 10$\%$ is needed to distinguish a hub from a non-hub in scale-free networks, generated for Simulation (iv). 

It is worth noting that if it is known that the true precision matrix is (approximately) block-diagonal, as in Simulation (iii), then computational speed-ups can be achieved by applying the proposed procedure \verb+hw.glasso+ to each block separately. In practice, as in \citep{Witten2011}, one may use a screening method on the elements of the sample covariance matrix $S_n$ to identify whether the solution to the \verb+hw.glasso+ problem (\ref{eq:HWGL1}) 
will be block-diagonal in which case the proposed method can be applied to each block separately.  

The simulation settings are considered for sample size $n=100$ with dimensions $p=50, 100, 200$, and sample size $n=250$ with dimensions $p=500, 1000$.

\subsection{A Two-Step Hubs Weighted Graphical Lasso \label{sec:TwoStepApproach}}

In our simulation studies, we observe that finite sample performance of the proposed \verb+hw.glasso+ can be improved by first identifying a set of candidate hubs $\widehat{\mathcal{H}}$ based on the \verb+hw.glasso+ estimate $\widehat{\Theta}_n$ 
and then penalizing the hub edges separately from the non-hub edges through a second weighted graphical lasso. In what follows, we outline this 2-step \verb+hw.glasso+ approach.

Based on the \verb+hw.glasso+ estimate $\widehat{\Theta}_n$, defined in (\ref{eq:HWGL1}), we identify a set of candidate hubs $\widehat{\mathcal{H}}$ (see the Remarks below for the choice of $\widehat{\mathcal{H}}$). We then construct a symmetric weight matrix $\widehat{W}=(\hat{w}_{ij})$, where 
\begin{align} \hat{w}_{ij} =  \begin{cases} 
      \lambda_1 & \text{if } i \in \widehat{\mathcal{H}} \text{ or } j \in \widehat{\mathcal{H}}, i \neq j  \\
      \lambda_2 & \text{if } i,j \notin \widehat{\mathcal{H}}, i \neq j  \\
      0 & \text{if }  i=j  \\
      \end{cases} 
      \label{eq:twostepweights}
\end{align}
for some tuning parameters $\lambda_1, \lambda_2>0$, and solve the weighted lasso optimization problem
\begin{align}
\bar{\Theta}_n = \displaystyle \operatorname*{arg\,max}_{\Theta \succ 0} \left\{  \log{\det{\Theta}} - \text{tr}(S_n \Theta) - \| \widehat{W} * \Theta \|_{1}  \right\} \nonumber,
\end{align}
where we refer to $\bar{\Theta}_n$ as the 2-step \verb+hw.glasso+ estimator of $\Theta$. The tuning parameter $\lambda_1$ controls the number of edges connecting a hub node to any other node in the graph, while the tuning parameter $\lambda_2$ controls the number of edges connecting two non-hub nodes. In our simulation studies, $\lambda_1$ and $\lambda_2$ are chosen using BIC.

An alternative choice of (\ref{eq:twostepweights}) is to use the adaptive weights $\hat{w}_{ij}^*=\lambda_1/|\tilde{\theta}_{ij}|^{\gamma_1}$, if $i \in \widehat{\mathcal{H}}$ or $j \in \widehat{\mathcal{H}}$, $i \neq j$,  and $\hat{w}_{ij}^*=\lambda_2/|\tilde{\theta}_{ij}|^{\gamma_1}$ if $i, j  \notin \widehat{\mathcal{H}}$, $i \neq j$, and 0 otherwise. More discussion is provided in Section 5.3.  \\

\noindent {\bf Remarks:}
\begin{enumerate}
\item[] Here we discuss two possible approaches for identifying a set of candidate hubs $\widehat{\mathcal{H}}$ based on the one-step \verb+hw.glasso+ estimate $\widehat{\Theta}_n$:
\begin{enumerate}
\item The set $\widehat{\mathcal{H}}$ can be obtained by setting a cutoff threshold for a node to be a hub. For example, as mentioned in Section \ref{sec:PerformanceMeasuresSettings}, we classify a node as a hub if it is connected to more than 10$\%$ of all other nodes. 
\item The set $\widehat{\mathcal{H}}$ can also be obtained by using a clustering approach. 
 From the first-step estimate $\widehat{\Theta}_n$, the degree of each node is computed and K-means clustering is then applied to cluster the nodes into two groups, where the hub group is characterized as the group with the larger mean degree. A similar approach based on a two-component Gaussian mixture model was considered by \cite{Charbonnier2010} in order to cluster nodes in a \emph{directed} graph as hubs and leaves. 
\end{enumerate}
\item[] In our simulation studies, we also consider the case where the hubs are known and thus take $\widehat{\mathcal{H}}=\mathcal{H}$.
\end{enumerate}

\subsection{Discussion of Simulation Results}

From Tables \ref{tab:Simulation1} to \ref{tab:Simulation3} corresponding to Simulations (i) to (iii), respectively, 
we see that when the true underlying network has hubs, the one-step \verb+hw.glasso+ procedure results in substantially better finite-sample performance compared to \verb+glasso+ and \verb+Ada-glasso+ that do not explicitly take hub structure into account. The \verb+hw.glasso+ procedure also outperforms the \verb+HGL+ 
and \verb+SF+ 
which are methods designed specifically for modelling networks with hubs. For Simulations (ii) and (iii) in which the hubs are not as highly connected, \verb+hw.glasso+  and \verb+Ada-glasso+ perform similarly when $n > p$, but the performance of \verb+hw.glasso+ increasingly improves relative to \verb+Ada-glasso+ as $p$ increases. The \verb+SF+ approach does not result in significant improvements over the \verb+glasso+ and \verb+Ada-glasso+ procedures, which is expected as it is not designed for estimating networks with very densely connected hubs. The \verb+HGL+ tends to perform better than \verb+glasso+ and \verb+Ada-glasso+ in terms of hub edge identification. However, with $c=0.5$ and $c=0.75$, which is a user-specified tuning parameter in the BIC-type quantity of \cite{Tan2014}, \verb+HGL+ leads to much denser graphs compared to \verb+glasso+. The tuning parameter $c$ controls the number of hubs in the graph, favouring more hubs when $c$ is small. Note that the value $c = 0.2$ is used by the authors in \cite{Tan2014}, which would have resulted in denser graphs and hence worse performance than what is shown here with $c=0.5$ or $c=0.75$. In practice, specifying an appropriate value of $c$ may require some prior knowledge of the number of hubs. 

Results of the 2-step \verb+hw.glasso+ procedures are also provided in Tables \ref{tab:Simulation1} to \ref{tab:Simulation3}. We see that when the true hubs are known in advance, which is a reasonable assumption in some biological applications, the 2-step \verb+hw.glasso+ that takes into account knowledge of these hubs results in significant improvements over the competing methods listed in the tables. While the 2-step 
\verb+hw.glasso+ in the case where the hubs are unknown performs better in higher dimensions than the one-step version, the 2-step procedure requires setting a cutoff threshold for a node to be considered a hub. Simulations were also conducted for the 2-step procedure using the adaptive weights $\hat{w}_{ij}^*$ introduced in Section 5.2, which yielded only slightly better results across all performance measures, and hence are not reported here.

In cases where $n \geq p$, we observe that the \verb+hw.glasso+ procedures (one-step and two-step methods) are better able to identify hub edges, leading to higher true positive rates compared to competing methods. In cases where $n<p$, all methods considered have greater difficulty in terms of edge identification. For Simulations (i) to (iii), the case $n=250$ and $p=500$, in particular, is challenging for all methods considered. Even in the case of the 2-step \verb+hw.glasso+, where the hubs are known in advance and less penalization is applied to hub edges, the true positive rate is low.

The results for Simulation (iv) are given in Table \ref{tab:Simulation4}. Note that 
the scale-free networks generated in this simulation have hubs that are not as highly connected as those in Simulations (i) to (iii). From Table \ref{tab:Simulation4}, it is thus not surprising that \verb+Ada-glasso+ performs well. When $p \geq n$, knowing the true hubs in advance and allowing for different levels of penalization between hub and non-hub edges, as in the 2-step \verb+hw.glasso+, results in better performance compared to the other methods across almost all performance measures. The one-step \verb+hw.glasso+ procedure performs well in terms of hub edge identification. The results for \verb+HGL+ are omitted as their method is not designed for estimating scale-free networks.  

To demonstrate the effect of $\gamma_1$ and $\gamma_2$ on the finite sample performance of \verb+hw.glasso+, we ran additional simulations and the results are summarized in Tables S1 and S2 of the Supplement. Table S1 corresponds to Simulations (i) and (ii), which cover case (a) of Proposition 1. Table S2 covers case (b) of the proposition and the simulation setting is described in the Supplement. In both tables, we first fix $\gamma_1=1$ and change $\gamma_2 \in \left\{0, 0.2, 0.5, 0.8, 1, 1.5, 2 \right\}$.  As expected, for smaller values of $\gamma_2$, the performance of \verb+hw.glasso+ is similar to that of \verb+Ada-glasso+. As $\gamma_2$ increases, \verb+hw.glasso+ outperforms \verb+Ada-glasso+ based on all the performance measures considered. On the other hand, as $\gamma_2$ increases beyond 1, the difference in performance by \verb+hw.glasso+ is minimal. This reaffirms our choice of $\gamma_1=\gamma_2=1$ in our simulations. In Table S2, we also consider the case where $\gamma_1=0$ and $\gamma_2 \in \left\{0.5, 0.8, 1, 1.5, 2 \right\}$. For this particular setting, the method performs comparably to the case where $\gamma_1>0$. In practice, we recommend using $\gamma_1>0$ and $\gamma_2>0$ in our proposed approach.

\section{Real data example \label{sec:MicrobiomeApplication}}

In this section, we illustrate the proposed methodology by estimating microbial interaction networks  
using undirected graphical models. The analysis is based on 
saliva microbiome relative abundance data sets of two \emph{Pan} species found in \cite{Li2013}.
We use relative abundances of genera in the saliva microbiomes of $n=23$ bonobos \emph{(Pan paniscus)} from the Lola ya Bonobo Sanctuary in the Democratic Republic of the Congo, and $n=22$ chimpanzees \emph{(Pan troglodytes)} from the Tacugama Chimpanzee Sanctuary in Sierra Leone. 

For the bonobos, 69 genera were identified along with 2 unknown/unclassified genera $(p=71)$. \emph{Enterobacter} (20.8$\%$) was the most abundant genus identified, followed by \emph{Porphyromonas} (10.3$\%$) and \emph{Neisseria} (9.7$\%$). For the chimpanzees, 79 genera were identified along with 2 unknown/unclassified genera $(p=81)$. The most abundant genera identified were \emph{Porphyromonas} (16.9$\%$), \emph{Fusobacterium} (14.0$\%$), \emph{Haemophilus} (11.4$\%$) and \emph{Neisseria} (8.1$\%$).

As microbial relative abundance data are compositional, after replacing zero abundance counts by 0.5, 
we use a centered log-ratio transformation \citep{Aitchison1981} of the data for our analysis. 
We then estimate undirected graphical models for each data set, using \verb+HGL+, \verb+Ada-glasso+, and 
 \verb+hw.glasso+ procedures. We also attempted \verb+SF+, but due to the small sample size $n$ relative to dimension $p$, this method had convergence issues and we did not obtain stable results, and thus it is not included here.
 For \verb+HGL+, we set $c=0.75$ and select its remaining three tuning parameters from fine grids. 
 For \verb+hw.glasso+, we select the tuning parameter 
$\lambda$ from a fine grid 
using BIC, and set $\gamma_1 = \gamma_2 = 1$. 
To obtain a graph that is reproducible under random sampling, we generate 
100 bootstrap samples and repeat the \verb+hw.glasso+ procedure on each sample. 
The stability of the network is then measured by the average proportion of edges 
reproduced by each bootstrap replicate. Only the edges that are reproduced in 
at least 80$\%$ of the bootstrap replicates are retained in the final network. 

Assuming hub structures for both the bonobo and chimpanzee microbial interaction networks and applying the \verb+hw.glasso+ procedure (retaining only reproducible edges), we found nodes corresponding to genera \emph{Actinobacillus}, \emph{Enterobacter} and \emph{Escherichia} to be highly connected for the bonobo group, and the node corresponding to genus \emph{Prevotella} to be highly connected for the chimpanzee group. There are 58 edges  that are common between the two groups. For both groups, there is a tendency for genera to correlate positively with other genera from the same phylum, especially within \emph{Proteobacteria} and \emph{Firmicutes}, which was also found in \cite{Li2013}. 

For each network, we use the \verb+R+ package \verb+igraph+ to evaluate several network measures, including network density, global clustering coefficient, betweenness centrality, and average path length. Differences in network measures between the bonobo and chimpanzee groups are assessed for statistical significance by permutation tests with 1000 randomizations. More specifically, we randomly assign the apes to one of two groups 1000 times. For each permutation, a network is estimated for each group and distributions of the differences in network indices are generated for statistical inference. No significant differences were found in terms of the global network structure (measured by global clustering coefficient, betweenness centrality and average path length) between the two groups. Significant differences in degree centrality were found for nodes corresponding to genera \emph{Escherichia} (0.35 v 0.14, $\text{p-value}=0.03$) and \emph{Peptostreptococcaceae} (0.02 v 0.18, $\text{p-value}=0.04$). 

\begin{figure}[tp]
\begin{center}
\includegraphics[scale=0.8]{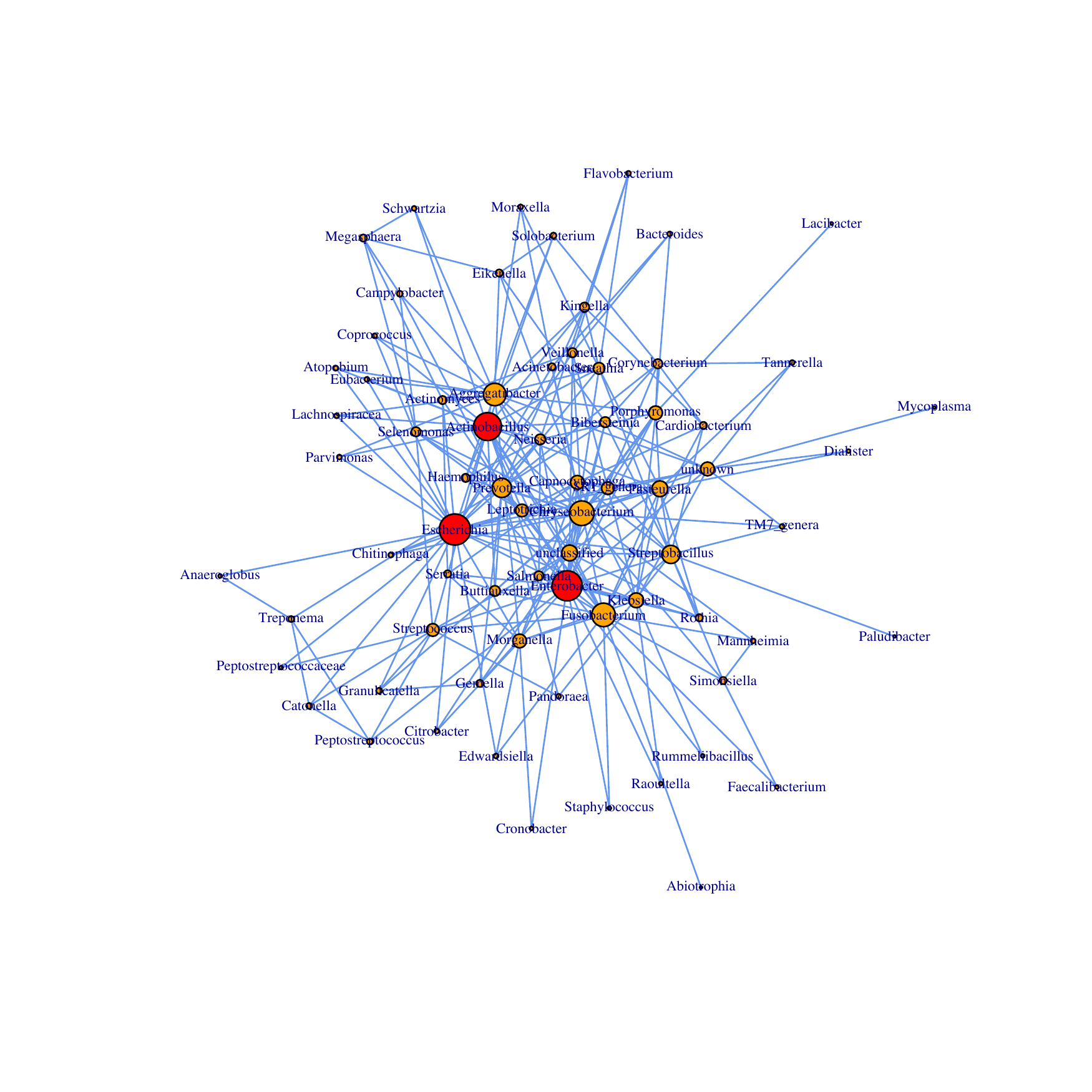}
\caption{Bonobo microbial interaction network, estimated by the hubs weighted graphical lasso.}
\label{fig:BonoboNetwork}
\end{center}
\end{figure}

\begin{figure}[tp]
\begin{center}
\includegraphics[scale=0.8]{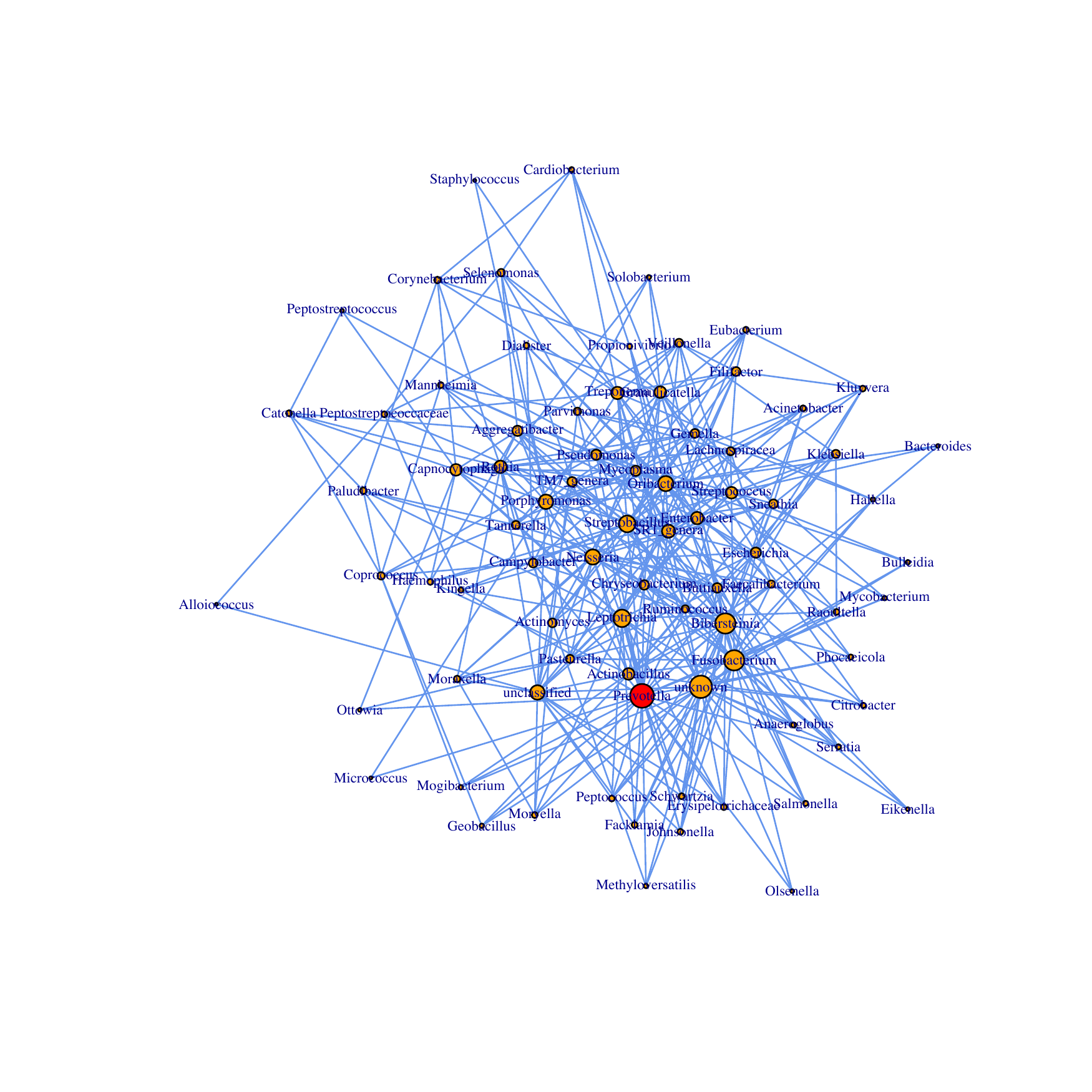}
\caption{Chimpanzee microbial interaction network, estimated by the hubs weighted graphical lasso.}
\label{fig:ChimpNetwork}
\end{center}
\end{figure}

The networks produced by our proposed method are displayed in Figures \ref{fig:BonoboNetwork} and \ref{fig:ChimpNetwork}. The hubs identified by our method were found to be the highest connected 
nodes by the \verb+Ada-glasso+, 
but our procedure assigned more edges to hubs and fewer edges to non-hubs, compared to the 
\verb+Ada-glasso+. 
The networks produced by \verb+HGL+ are displayed in Figures S1 and S2 of the Supplement. \verb+HGL+ also identified \emph{Actinobacillus} and \emph{Enterobacter} as highly connected nodes for the bonobo group, along with genera \emph{Acinetobacter}, \emph{Streptobacillus}, \emph{Sneathia}, \emph{Aggregatibacter} and \emph{Bacteroides}. For the chimp group, \verb+HGL+ identified \emph{Prevotella} as the highest connected node along with  \emph{Anaeroglobus}, \emph{Ruminococcus}, \emph{Faecalibacterium}, and followed by \emph{Salmonella} and \emph{Sneathia}. This tendency of \verb+HGL+ to identify more nodes as hubs compared to \verb+hw.glasso+ and other competing methods was also observed in our simulation studies.

\section{Conclusions \label{sec:Conlusions}}


In this paper, we proposed a new weighted graphical lasso approach for estimating networks with hubs that makes use of informative weights that allow for hub structure. We showed that the proposed method, referred to as the hubs weighted graphical lasso (\verb+hw.glasso+), is both estimation and selection consistent. We then demonstrated with simulated data that the proposed method performs significantly better than methods that do not explicitly take hub structure into account, but it also outperforms network estimation procedures designed for modelling networks with hubs, such as the \verb+HGL+ of \cite{Tan2014} and the re-weighted $L_1$-regularization approach of \cite{LiuIhler2011}. The former is designed for estimating networks with very densely connected hub nodes, referred to as ``super hubs'', while the latter is designed for estimating scale-free networks, for which there may be no clear distinction between \emph{hub} and \emph{non-hub} nodes. Our proposed method can accommodate both networks with so-called ``super hubs'' as well as scale-free networks. 

In the proposed method, the construction of the weights in (\ref{eq:Weights}) requires an initial consistent estimator of the precision matrix $\Theta_0$. Under the regularity conditions in Theorem 1, the standard \verb+glasso+ estimator provides such an estimator. When these conditions are violated or an initial consistent estimator is not available, properties of our proposed estimator as well as the roles of the tuning parameters $\gamma_1$ and $\gamma_2$ in the weights are presently unknown to us. Such cases are the subject of future research. Another possible research direction is to investigate the extension of our method to group selection \cite{Bach2008}, which will depend on the definition of grouping in the context of estimating a sparse precision matrix and may require a re-design of the penalty function. 

Our current work focuses on the problem of static network modelling, where the inferred  
network may provide a snapshot of a network structure at a single time point. 
In some applications, networks may undergo changes over time in 
response to changes in external conditions  
and the temporal variation of these networks can be captured by dynamic networks \citep{Faust2015}. 
Techniques developed for static network modelling will pave the way for proposing new 
approaches for modelling the dynamics of networks with hubs.

\section*{Acknowledgments}
\noindent
The authors would like to thank 
the Natural Sciences and Engineering Research Council of Canada and 
the Fonds de recherche du Qu\'ebec-Nature et technologies for partly supporting this research. 

\section*{Appendix \label{sec:appendix}}

Let $A^+ = \text{diag}(A)$ be a diagonal matrix with the diagonal elements of a matrix $A$, 
and further let $A^- = A - A^+$. Also, for any two matrices $B_{m \times m}  =\{b_{ij}\}$ and 
$C_{n \times n} = \{c_{kl}\}$, their Kronecker product is $B \otimes C  = \{b_{ij} c_{kl} \}$.
For any closed bounded convex set $\mathcal{C}$ 
which contains ${\bf 0}$, its boundary is denoted by $\partial {\cal C}$.
Recall from (\ref{eq:Weights}) that 
$\|\boldsymbol{\theta}^{0}_{\neg i}  \|_1=\sum_{k \neq i} |\theta_{ik}^{0}|$. 
We also use the result of Lemma 3 of \cite{BicLev2008}, re-stated in what follows.

\begin{lemma}
\label{lem1}
Let ${\bf X}_i = (X_{i1}, X_{i2}, \ldots, X_{ip_n})^\top, i= 1, 2, \ldots, n$, be iid Gaussian with mean ${\bf 0}$ and covariance matrix $\bSig_0 = \{\sigma^0_{ij}\}$ such that 
$\phi_{\text{max}}(\bSig_0) \le \tau^{-1}_1 < \infty$, then 
\[
P\bigg [ \bigg | \sum_{i=1}^n (X_{ij} X_{ik} - \sigma^0_{ij})\bigg | \ge n \nu \bigg ] \le c_1 \exp(-c_2 n \nu^2)~~,~~
|\nu| \le \delta,
\]
where $c_1, c_2$ and $\delta$ depend on $\tau_1$ only. 
\end{lemma}

We now proceed with the proof of our first result. \\

\noindent {\bf Proof of Theorem \ref{thm:consistency}}:
The idea of the proof is inspired by the proof of Theorem 1 of \cite{Rothman2008}. Here we work with 
the negative penalized log-likelihood $Q_n(\Theta) = -pl_n(\Theta)$. 
Let $\Delta = \Theta - \Theta_0$, and define $G(\Delta) = Q_n(\Theta_0 + \Delta) - Q_n(\Theta_0)$
which is a convex function of $\Delta$. 
Also let $\widehat{\Delta}_n = \widehat \Theta_n - \Theta_0$, 
where $\widehat \Theta_n$ is the weighted \verb+glasso+ estimator 
which minimizes $Q_n(\Theta)$ or equivalently $\widehat{\Delta}_n$
minimizes $G(\Delta)$. 
Then $G(\widehat{\Delta}_n) \leq G({\bf 0}) = 0$. 
Now if we take a closed bounded convex set $\mathcal{C}$ 
which contains ${\bf 0}$, and show that $G$ is strictly positive everywhere on the 
boundary $\partial \mathcal{C}$, then it implies that $G$ has 
its minimizer $\widehat \Delta_n$ inside $\mathcal{C}$ since $G$ is continuous 
and $G({\bf 0})=0$. 
Define the set
\begin{align*}
\mathcal{C} = \left\{ \Delta:  \Delta = \Delta^{\top}, \| \Delta \|_F^2 \leq M r_n^2  \right\} 
\end{align*}
with the boundary 
\[
\partial \mathcal{C} = \left\{ \Delta:  \Delta = \Delta^{\top}, \| \Delta \|_F^2 = M r_n^2  \right\},
\]
where $M$ is a positive constant and $r_n = \left\{ (p_n+q_n) (\log{p_n})/n \right\}^{1/2}$. Then we must show that $P( \displaystyle \inf_{\Delta \in \partial \mathcal{C}} G(\Delta) > 0 ) \rightarrow 1$, as $n \rightarrow \infty$. 
We proceed as follows. 

Using \eqref{eq:like} and \eqref{eq:penlike}, we have that
\begin{align}
G(\Delta) = Q_n(\Theta_0 + \Delta) - Q_n(\Theta_0) & = 
-\log{\det{(\Theta_0 + \Delta)}} + \log{\det{(\Theta_0)}} + \text{tr}(S_n (\Theta_0+\Delta)) - \text{tr}(S_n \Theta_0)  \nonumber \\
&+ \lambda_n \left\{ \| W*(\Theta_0 +  \Delta) \|_1 - \| W*\Theta_0 \|_1 \right\}. \nonumber
\end{align}

Now, using the Taylor expansion of $f(t) = \log{\det{(\Theta + t \Delta)}}$ and the fact that 
$\Delta, \bSig_{0}$ and $\Theta_0$ are all symmetric matrices, we have 
\small
\[
\log{\det{(\Theta_0 +  \Delta)}} - \log{\det{(\Theta_0)}} 
= 
\text{tr}(\bSig_0 \Delta) - [\text{vec}({\Delta})]^\top \left\{ \int_{0}^1 (1 - v)(\Theta_0 + v \Delta)^{-1} \otimes (\Theta_0 + v \Delta)^{-1}  \,dv \right\} \text{vec}(\Delta),
\]
where $\text{vec}({\Delta})$ is the vectorized version of the matrix $\Delta$ to match the multiplication. 
Thus,  
\begin{align*}
G(\Delta) & =   \text{tr}(\Delta(S_n - \bSig_0)) +  [\text{vec}({\Delta})]^\top \left\{ \int_{0}^1 (1 - v)(\Theta_0 + v \Delta)^{-1} \otimes (\Theta_0 + v \Delta)^{-1}  \,dv \right\} \text{vec}(\Delta)  \\ 
&+ \lambda_n \left\{ (\| W*( \Theta_0 +  \Delta) \|_1 - \| W*\Theta_0 \|_1 \right\} \nonumber 
= {\bf I_1} + {\bf I_2} + {\bf I_3}. 
\end{align*}
\normalsize

To show that $G({\Delta})$ is strictly positive on $\partial \mathcal{C}$, we need to bound 
the quantities ${\bf I_1}$, ${\bf I_2}$ and ${\bf I_3}$.  
By the union sum inequality, the Cauchy-Schwartz inequality, and Lemma \ref{lem1}, 
there exist positive constants $C_1$ and $C_2$ such that with probability tending to 1, 
as $n \to \infty$,
\begin{align}
- | {\bf I_1} |  = 
 - \Big | \text{tr}(\Delta(S_n - \bSig_0)) \Big | &\geq 
 - \Big|  \sum_{i \neq j}^{p_n} \Delta_{ij} (s_{ij} - \sigma^0_{ij}) \Big| 
 - \Big|  \sum_{i=1}^{p_n} \Delta_{ii} (s_{ii} - \sigma^0_{ii})   \Big| \nonumber \\
& \geq  - \max_{i \not= j} |s_{ij} - \sigma^0_{ij}| \times \| \Delta^{-} \|_1 - 
\sqrt{p_n} \max_{1\le i \le p_n} |s_{ii} - \sigma^0_{ii}| \times \| \Delta^{+} \|_F
\nonumber \\
& \geq  - C_1 \left( \frac{\log{p_n}}{n} \right)^{1/2} \|  \Delta^{-} \|_1 - C_2 \left( \frac{p_n\log{p_n}}{n}\right)^{1/2} \|  \Delta^{+} \|_F  \nonumber \\
& \geq 
- C_1 \left( \frac{\log{p_n}}{n} \right)^{1/2} \left( \| {\Delta}_T^{-} \|_1 
- \| {\Delta}_{T^c}^{-} \|_1 \right) - C_2 \left( \frac{(p_n+q_n)\log{p_n}}{n}\right)^{1/2} \|  \Delta^{+} \|_F.
\label{ineq-I1}
\end{align}

Also, as in \cite{Rothman2008}, we have that 
\begin{align*}
& \phi_{\text{min}} \left\{ \int_{0}^1 (1 - v)(\Theta_0 + v \Delta)^{-1} \otimes (\Theta_0 + v \Delta)^{-1}  \,dv \right\} 
 \ge 
\int_{0}^1 (1 - v) \phi^2_{\text{min}}\{ (\Theta_0 + v \Delta)^{-1} \} dv \\ 
& \ge 
\frac{1}{2} \min_{0 \le v \le 1} \phi^2_{\text{min}}\{ (\Theta_0 + v \Delta)^{-1} \}
 \ge \frac{1}{2} \min \bigg \{\phi^2_{\text{min}}\{ (\Theta_0 + \Delta)^{-1}\}: \Delta \in \mathcal{C} \bigg \} \\
& = 
\frac{1}{2} \min \bigg \{\phi^{-2}_{\text{max}}\{ (\Theta_0 + \Delta)\}: \Delta \in \mathcal{C} \bigg \} 
 \geq \frac{1}{4\tau_2^2}
\end{align*}
with probability tending to 1, where the last inequality is due to the regularity Condition 1, and 
also the fact that for all $\|\Delta\| \in \mathcal{C}$, we have $\|\Delta\| \le \|\Delta\| \le M r_n^2 = o(1)$, 
as $n \to \infty$. Thus, using the above inequality, with probability tending to 1, we have 
\begin{align}
 \label{ineq-I2}
|{\bf I_2}| = \bigg | 
[\text{vec}(\Delta)]^\top \left\{ \int_{0}^1 (1 - v)(\Theta_0 + v \Delta)^{-1} \otimes (\Theta_0 + v \Delta)^{-1}  \,dv \right\} \text{vec}(\Delta) \bigg | \geq \frac{1}{4\tau_2^2} \| \Delta \|_F^2.
\end{align}

Next, since the penalty is decomposable \citep{Negahban2012}, we have that
\begin{align}
\label{ineq-I3}
| {\bf I_3} | \geq \lambda_n \left( \| W* {\Delta}_{T^c}^{-} \|_1 - \| W* {\Delta}_T^{-} \|_1 \right).
\end{align}

Therefore, using the fact that $\| \Delta^-_T \|_1 \leq \sqrt{q_n}~ \| \Delta^- \|_F \leq \sqrt{p_n+q_n}~ \| \Delta^- \|_F$, and \eqref{ineq-I1}-\eqref{ineq-I3}, we find that for large $n$,
\begin{align}
G(\Delta) 
& \geq \frac{\gamma}{4\tau_2^2} \| \Delta  \|_F^2 -  C_2  \left\{ \frac{ (p_n+q_n) \log{p_n}}{n}\right\}^{1/2} \| \Delta^{+} \|_F  +  \left\{  \lambda_n  \min_{(i,j) \in T^c} \widetilde{w}_{ij} - C_1 \left( \frac{ \log{p_n}}{n}\right)^{1/2} \right\} \| \Delta_{T^c}^{-} \|_1\nonumber \\
&  - \left\{  \lambda_n \max_{(i,j) \in T} \widetilde{w}_{ij} + C_1 \left( \frac{ \log{p_n}}{n}\right)^{1/2} \right\} \| \Delta_{T}^{-} \|_1  \nonumber \\ 
& \geq \frac{\gamma}{4\tau_2^2} \| \Delta  \|_F^2 -  C_2  \left\{ \frac{ (p_n+q_n) \log{p_n}}{n}\right\}^{1/2} \| \Delta^{+} \|_F  - \left\{  \lambda_n \max_{(i,j) \in T} \widetilde{w}_{ij} + C_1 \left( \frac{ \log{p_n}}{n}\right)^{1/2} \right\} \| \Delta_{T}^{-} \|_1  \nonumber 
\end{align}
due to the condition 
$\left\{ (\log{p_n})/n \right\}^{1/2} \left( \min_{(i,j) \in T^c} \widetilde{w}_{ij} \right)^{-1}=O_p(\lambda_n)$. \\

Now, using $\Delta = \Delta^+ + \Delta^-$, we have that
\begin{align}
G(\Delta)  &\geq \left\{ \frac{\gamma}{4\tau_2^2}  -  \Bigg[  \lambda_n \max_{(i,j) \in T} \widetilde{w}_{ij} + 
C_1 \left( \frac{\log{p_n}}{n} \right)^{1/2}  \Bigg] (p_n+q_n)^{1/2}    \| \Delta^-  \|_F^{-1} \right\}   \| \Delta^-  \|_F^2 \nonumber \\
& + \left\{  \frac{\gamma}{4\tau_2^2} - C_2\left( \frac{(p_n+q_n)\log{p_n}}{n} \right)^{1/2} \| \Delta^+ \|_F^{-1}  \right\} \| \Delta^+ \|_F^2 \nonumber \\
& = \left\{ \frac{\gamma}{4\tau_2^2}  -  \Bigg[  \frac{\lambda_n \max_{(i,j) \in T} \widetilde{w}_{ij} }{C_1 (\log{p_n}/n)^{1/2}  }   + 1 \Bigg] C_1\left( \frac{(p_n+q_n) \log{p_n}}{n} \right)^{1/2}   \| \Delta^-  \|_F^{-1} \right\}   \| \Delta^-  \|_F^2 \nonumber \\
& + \left\{  \frac{\gamma}{4\tau_2^2} - C_2\left( \frac{(p_n+q_n)\log{p_n}}{n} \right)^{1/2} \| \Delta^+ \|_F^{-1}  \right\} \| \Delta^+ \|_F^2. \nonumber
\end{align} 

Since $\lambda_n \max_{(i,j) \in T} \widetilde{w}_{ij} = O_p\left( \left\{ (\log{p_n})/n \right\}^{1/2}  \right)$, 
then on the boundary set $\partial \mathcal{C}$, where $\| \Delta  \|_F^2 =  M r_n^2$, 
we have that
\begin{equation}
\begin{split}
G(\Delta) 
&\geq \|  \Delta^-  \|_F^2 \left( \frac{\gamma}{4\tau_2^2} -  C_1/\sqrt{M} \right) + \|  \Delta^+  \|_F^2 \left( \frac{\gamma}{4\tau_2^2} -  C_2/\sqrt{M} \right). \nonumber
\end{split}
\end{equation} 
Thus, for $M$ sufficiently large, we have that $G(\Delta)>0$ for any 
$\Delta \in \partial \mathcal{C}$, which completes the proof. 
$\blacksquare$\\


The result of the following Lemma is used for proving Theorem \ref{thm:sparsistency}. 
First, we introduce some notation. We write the true precision matrix $\Theta_0$ as a $p_n(p_n+1)/2$-dimensional 
vector $\bpsi_0$ by taking $\bpsi_0 := \bpsi_0(\Theta_0) = (\bpsi_{01}, \bpsi_{02})$ such that $\bpsi_{02}={\bf 0}$. 
A similar presentation is used for any precision matrix $\Theta$: $\bpsi = \bpsi(\Theta)=(\bpsi_1, \bpsi_2)$. 
Recall $r_n = \left\{ (p_n+q_n) (\log{p_n})/n \right\}^{1/2}$.

\begin{lemma}
\label{lem2}
For any precision matrix $\Theta$ such that 
\begin{align}
\| \Theta - \Theta_0 \|_{\text F} &= O_p(r_n)
~~,~~\| \Theta - \Theta_0 \|^2 = O_p(\eta_n)
 \nonumber 
\end{align}
where $\eta_n \rightarrow 0$ as $n \rightarrow \infty$, if $\left\{ \frac{\log{p_n}}{n} + \eta_n \right\} \left\{ \min_{(i,j) \in T^c} \widetilde{w}_{ij}\right\}^{-2} = O_P(\lambda_n^2)$, then 
\begin{align}
pl_n((\bpsi_1, \bpsi_2); \lambda_n) - pl_n((\bpsi_1, {\bf 0}); \lambda_n) < 0
\end{align}
with probability tending to 1 as $n \rightarrow \infty$, where $\bpsi = \bpsi(\Theta)=(\bpsi_1, \bpsi_2)$.  
\end{lemma}

\noindent 
{\bf Proof of Lemma \ref{lem2}:} 
Recall the definition of the penalized log-likelihood in \eqref{eq:penlike} and with the general weighted $L_1$-penalty 
in \eqref{eq:HWGL1}. We have that
\begin{align}
\label{diff:plike}
pl_n((\bpsi_1, \bpsi_2); \lambda_n) - pl_n((\bpsi_1, {\bf 0}); \lambda_n) & = \left\{ \ell_n((\bpsi_1, \bpsi_2)) - \ell_n((\bpsi_1, {\bf 0})) \right\} - \sum_{(i,j) \in T^c} \lambda_n \widetilde{w}_{ij} | \theta_{ij} |.
\end{align}

We first analyze the difference in the log-likelihood part. By the Mean Value Theorem,
\begin{align}
\label{diff:like}
\ell_n((\bpsi_1, \bpsi_2)) - \ell_n((\bpsi_1, {\bf 0})) & = \Big[\frac{\partial \ell_n((\bpsi_1, \bxi))}{\partial \bpsi_2} \Big]^{\top} \times \bpsi_2,
\end{align}
where $\bxi$ is a vector between $\bpsi_2$ and $\bpsi_{02}={\bf 0}$ such that 
$\| \bxi \|_2 \leq \| \bpsi_2 \|_2$.

As in the proof of Theorem 2 of \cite{LamFan2009}, we have that
\begin{align*}
\frac{\partial \ell_n((\bpsi_1, \bxi))}{\partial \bpsi_2} & = \left\{ \sigma_{ij}(\bxi) - s_{ij}: (i,j) \in T^c \right\}.
\end{align*}
We need to assess the orders of $\sigma_{ij}(\bxi) - s_{ij}$ as $n \rightarrow \infty$. Note that 
\begin{align*}
\sigma_{ij}(\bxi) - s_{ij} = (\sigma_{ij}(\bxi) - \sigma_{ij}^0) + (\sigma_{ij}^0 - s_{ij}).
\end{align*}
By \cite{LamFan2009}, $| \sigma_{ij} - \sigma_{ij}^0 | \leq \| \bSig - \bSig_0 \|$, which has the order
\begin{align}
\| \bSig - \bSig_0 \| = \| \bSig (\Theta - \Theta_0) \bSig_0 \| \leq \| \bSig \| \times \| \Theta - \Theta_0 \| \times \| 
\bSig_0 \| \nonumber
\end{align}
and $\| \bSig_0 \| = O(1)$ by Condition 1. Also using that $\eta_n \rightarrow 0$ so that $\lambda_{\text{min}}(\Theta - \Theta_0) = o(1)$ for $\|\Theta - \Theta_0 \| = O(\eta_n^{1/2})$, we find
\begin{align*}
\| \bSig \| = \lambda_{\text{min}}^{-1}(\Theta) \leq \Big[ \lambda_{\text{min}}(\Theta_0) + \lambda_{\text{min}}(\Theta-\Theta_0) \Big]^{-1} = O(1)
\end{align*} 
and since $\| \Theta - \Theta_0 \| = O(\eta_n^{1/2})$, 
we have that $| \sigma_{ij} - \sigma_{ij}^0 | = O(\eta_n^{1/2})$.

Since $\sigma_{ij}(\bxi)$ is between $\sigma_{ij}$ and $\sigma_{ij}^0$, $| \sigma_{ij}(\bxi) - \sigma_{ij}^0 | = O(\eta_n^{1/2})$. Therefore,
\begin{align}
\max_{(i,j) \in T^c} | \sigma_{ij}(\bxi) - s_{ij} | = O_p( |s_{ij} - \sigma_{ij}^0 | + \eta_n^{1/2}) 
\label{eq:1}
\end{align}
as $n \rightarrow \infty$. On the other hand, by Lemma \ref{lem1}, 
\begin{align}
\max_{(i,j) \in T^c} | s_{ij} - \sigma_{ij}^0 | = O_p\left\{ \left(\frac{\log{p}}{n}\right)^{1/2} \right\}
\label{eq:2}
\end{align}
as $n \rightarrow \infty$. Equations \eqref{eq:1} and \eqref{eq:2} imply that, as $n \rightarrow \infty$, 
\begin{align}
\label{order1}
\max_{(i,j) \in T^c} | \sigma_{ij}(\bxi) - s_{ij} | 
= 
O_p\left\{ \left(\frac{\log{p}}{n}\right)^{1/2} + \eta_n^{1/2}
\right\}.
\end{align}

Going back to the log-likelihood difference in \eqref{diff:like}, it can be written as
\begin{align}
\label{diff2:like}
 \ell_n((\bpsi_1, \bpsi_2)) - \ell_n((\bpsi_1, {\bf 0}))  = \sum_{(i,j) \in T^c} \left\{ \sigma_{ij}(\bxi) - s_{ij} \right\}|\theta_{ij}|.
\end{align}
Replacing the order assessment \eqref{order1} in \eqref{diff2:like}, we have that 
\begin{align*}
pl_n((\bpsi_1, \bpsi_2); \lambda_n) - pl_n((\bpsi_1, {\bf 0}); \lambda_n) & = 
 \left\{ \ell_n((\bpsi_1, \bpsi_2)) - \ell_n((\bpsi_1, {\bf 0})) \right\} - \sum_{(i,j) \in T^c} \lambda_n \widetilde{w}_{ij} | \theta_{ij} |  \\ 
 & \leq \sum_{(i,j) \in T^c}  \left\{ O_p\left((\log{p}/n)^{1/2} + 
 \eta_n^{1/2} \right)|\theta_{ij} | - \lambda_n \tilde{w}_{ij} | \theta_{ij} |  \right\} \\ 
 & = \sum_{(i,j) \in T^c}  \left\{ O_p\left((\log{p}/n)^{1/2} + \eta_n^{1/2} \right)- \lambda_n \tilde{w}_{ij}  \right\} |\theta_{ij} |  < 0 
\end{align*}
if $\lambda_n^2 > \left(\min_{(i,j) \in T^c} \tilde{w}_{ij}\right)^{-2}\left( \frac{\log{p_n}}{n} + \eta_n \right)$, 
for large $n$. 
In other words, if $\left( \frac{\log{p_n}}{n} + \eta_n \right) \left\{ \min_{(i,j) \in T^c} \widetilde{w}_{ij}\right\}^{-2} = O_p(\lambda_n^2)$, then with probability approaching 1, as $n \rightarrow \infty$,
\begin{align*}
pl_n((\bpsi_1, \bpsi_2); \lambda_n) - pl_n((\bpsi_1, {\bf 0}); \lambda_n) <0
\end{align*}
and this completes the proof. $\blacksquare$\\

The implication of Lemma \ref{lem2} is that in the neighbourhood (specified by the conditions of this Lemma)  
of the true precision matrix, $\bpsi_0 = (\bpsi_{01}, \bpsi_{02}) = (\bpsi_{01}, {\bf 0})$, the penalized 
log-likelihood function $pl((\bpsi_1, \bpsi_2); \lambda_n)$ 
is maximized only when $\bpsi_2 = {\bf 0}$.
We now proceed to the proof of Theorem \ref{thm:sparsistency}. \\

\noindent {\bf Proof of Theorem \ref{thm:sparsistency}: } Let $(\widehat{\bpsi}_{n1},{\bf 0})$ 
be the maximizer of $pl_n((\bpsi_1,{\bf 0}); \lambda_n)$ which is considered 
as a function of $\bpsi_1$ only. Then in the neighbourhood
\begin{align*}
\| \Theta - \Theta_0 \|_F &= O_p\left\{ \left( \frac{(p_n+q_n)\log{p_n}}{n}  \right)^{1/2}  \right\}   \\
\| \Theta - \Theta_0 \|^2 &= O_p(\eta_n),
\end{align*}
we have that, by Lemma \ref{lem2},
\begin{align*}
pl_n((\bpsi_1, \bpsi_2); \lambda_n) - pl_n((\widehat{\bpsi}_{n1}, 0); \lambda_n) &= 
\{ pl_n((\bpsi_1, \bpsi_2); \lambda_n) - pl_n(\bpsi_{1}, {\bf 0}); \lambda_n) \}  \nonumber \\
&+ 
\{ pl_n((\bpsi_1, {\bf 0}); \lambda_n) - pl_n((\widehat{\bpsi}_{n1}, {\bf 0}); \lambda_n) \} < 0 
\end{align*}
with probability tending to 1 as $n \rightarrow \infty$. 
Therefore, in the chosen neighbourhood of $\Theta_0$, with probability tending to 1 as $n \rightarrow \infty$, 
the maximum of $pl_n((\bpsi_1, \bpsi_2); \lambda_n)$ indeed happens at $(\widehat{\bpsi}_{n1},{\bf 0})$.
$\blacksquare$ \\

\noindent 
{\bf Proof of Proposition \ref{prop1}}: 
Theorems \ref{thm:consistency} and \ref{thm:sparsistency} 
require choices of the tuning parameter $\lambda_n$ and the (possibly random) weights $\widetilde w_{ij}$  
 that, as $n \to \infty$, satisfy conditions
\begin{align}
\label{con1}
     \lambda_n \max_{(i,j) \in T} \widetilde{w}_{ij} & = O_p \left( \sqrt{\frac{\log{p_n}}{n}} \right) \\
     \label{con2}
    \left( \frac{\log{p_n}}{n} \right)^{1/2} \left\{  \min_{(i,j) \in T^c} \widetilde{w}_{ij} \right\}^{-1}  
    & = O_p(\lambda_n) \\
    \label{con3}
     \eta_n  \left\{  \min_{(i,j) \in T^c} \widetilde{w}_{ij} \right\}^{-2}  &= O_p(\lambda_n^2), 
\end{align}
where $\lambda_n , \eta_n \rightarrow 0$. We now verify these conditions for the suggested weights 
$\widetilde w_{ij}$ in \eqref{eq:Weights} used in the hubs weighted graphical lasso (\verb+hw.glasso+). 
Note that these weights are constructed based on the popular graphical lasso (\verb+glasso+) 
estimator $\widetilde{\Theta}_n$. By \cite{Rothman2008}, we have that as $n \to \infty$, 
\begin{align}
\label{rate1}
\| \widetilde{\Theta}_n - \Theta_0 \|_F = O_p\left(\sqrt{ \frac{(p_n+q_n)\log{p}_n}{n} } \right).
\end{align}

We start with \eqref{con1}. By the definitions of the weights in \eqref{eq:Weights}, we have 
\begin{align*}
\lambda_n \max_{(i,j) \in T} \widetilde{w}_{ij} &= 
  \max_{(i,j) \in T} \frac{\lambda_n}{ | \tilde{\theta}_{ij} |^{\gamma_1} \left\{  \| \tilde{\boldsymbol{\theta}}_{\neg i}  \|_1 \cdot   \| \tilde{\boldsymbol{\theta}}_{\neg j} \|_1 \right \}^{\gamma_2}} 
   \\
& =  \frac{\lambda_n}{\min_{(i,j) \in T}  |\tilde{\theta}_{ij}|^{\gamma_1} \left\{ \|  \tilde{\boldsymbol{\theta}}_{\neg i} \|_1 \| \tilde{\boldsymbol{\theta}}_{\neg j}  \|_1 \right\}^{\gamma_2}} \le \lambda_n C(\tau_3),
\end{align*}
where the last inequality is due to \eqref{rate1}, the regularity Condition 2 and that $T \neq \emptyset$, 
and $C(\tau_3)$ is a function of $\tau_3$. 
 Thus \eqref{con1} is satisfied, if we choose $\lambda_n$ as 
 \begin{align*}
\lambda_n = O\left( \sqrt{ \frac{\log{p_n}}{n} }\right),
\end{align*}
 as required in 
 \eqref{final-condition1a} and \eqref{final-condition1b} 
 of the two parts {\bf (a)-(b)} of the Proposition. 

For \eqref{con2}-\eqref{con3}, we divide the proof as follows. 

\noindent
{\bf Part (a)}. By the conditions of this part of the Proposition, 
since there exists a pair $(i,j) \in T^c$ such that 
$\|\boldsymbol{\theta}^{0}_{\neg i}  \|_1 \neq 0$ 
and $\|\boldsymbol{\theta}^{0}_{\neg j}  \|_1 \neq 0$, then by using \eqref{rate1}, Condition 2, 
and that $T \neq \emptyset$, we have that for large $n$,
\begin{align*}
 \left\{\min_{(i,j) \in T^c} \widetilde{w}_{ij}  \right\}^{-1} &= 
  \max_{(i,j) \in T^c} \left\{  | \tilde{\theta}_{ij}|^{\gamma_1} 
  \Big[ \| \tilde{\boldsymbol{\theta}}_{\neg i}  \|_1 \| \tilde{\boldsymbol{\theta}}_{\neg j}  \|_1 \Big]^{\gamma_2} \right\} \\
& = \xi_2(\tau_3) \max_{(i,j) \in T^c}  | \tilde{\theta}_{ij}|^{\gamma_1} 
\leq \xi_2(\tau_3) \left\{ \sum_{(i,j) \in T^c} \tilde{\theta}_{ij}^2 \right\}^{\gamma_1/2} \\
& \le 
\xi_2(\tau_3) \left\{ \frac{(p_n+q_n)\log{p_n}}{n} \right\}^{\gamma_1/2}
\end{align*}
for some constant $\xi_2(\tau_3)>0$. Thus, to satisfy \eqref{con2}-\eqref{con3}, 
and using the above inequality, it is sufficient to choose $\lambda_n$ and $\eta_n$ such that  
\begin{align*}
\left( \frac{\log{p_n}}{n} \right)^{1/2} \left\{ \frac{(p_n+q_n)\log{p_n}}{n} \right\}^{\gamma_1/2} \lambda^{-1}_n  & = O(1) \\
\sqrt \eta_n \left( \frac{(p_n+q_n)\log{p_n}}{n} \right)^{\gamma_1/2} \lambda_n^{-1} & = O(1)
\end{align*}
as required by \eqref{final-condition1a} and \eqref{final-condition2a} of the Proposition. \\

\noindent
{\bf Part (b)}. By the conditions of this part of the Proposition, since for all $(i,j) \in T^c$, either $\|\boldsymbol{\theta}^{0}_{\neg i}  \|_1 = 0$ or $\|\boldsymbol{\theta}^{0}_{\neg j}  \|_1 = 0$, we have that by \eqref{rate1} and that $T \neq \emptyset$,
for large $n$,  
\begin{align*}
 \left\{\min_{(i,j) \in T^c} \widetilde{w}_{ij}  \right\}^{-1} & = 
  \max_{(i,j) \in T^c} \left\{  | \tilde{\theta}_{ij}|^{\gamma_1} 
  \Big[ \| \tilde{\boldsymbol{\theta}}_{\neg i}  \|_1 \| \tilde{\boldsymbol{\theta}}_{\neg j}  
  \|_1 \Big]^{\gamma_2} \right\} \le \left\{ \sum_{(i,j) \in T^c} \tilde{\theta}_{ij}^2 \right\}^{\gamma_1/2} 
  \max_{(i,j) \in T^c} \left\{
  \Big[ \| \tilde{\boldsymbol{\theta}}_{\neg i}  \|_1 \| \tilde{\boldsymbol{\theta}}_{\neg j}  
  \|_1 \Big]^{\gamma_2} \right\}
  \\
& 
\leq \xi_2(\tau_3) \left\{ \frac{(p_n+q_n)\log{p_n}}{n} 
\right\}^{\gamma_1/2} 
\left\{ \frac{p_n (p_n+q_n)\log{p_n}}{n}
\right\}^{\gamma_2/2}
\end{align*}
for some constant $\xi_2(\tau_3)>0$. Thus, to satisfy \eqref{con2}-\eqref{con3}, 
and using the above inequality, it is sufficient to choose $\lambda_n$ and $\eta_n$ such that  
\begin{align*}
\left( \frac{\log{p_n}}{n} \right)^{1/2} \left\{ \frac{(p_n+q_n)\log{p_n}}{n} \right\}^{\gamma_1/2} 
\left\{ \frac{p_n (p_n+q_n)\log{p_n}}{n}
\right\}^{\gamma_2/2}
\lambda^{-1}_n  & = O(1) \\
\sqrt \eta_n \left( \frac{(p_n+q_n)\log{p_n}}{n} \right)^{\gamma_1/2} 
\left\{ \frac{p_n (p_n+q_n)\log{p_n}}{n}
\right\}^{\gamma_2/2}
\lambda_n^{-1} & = O(1)
\end{align*}
as required by \eqref{final-condition1b} and \eqref{final-condition2b} of the Proposition.

\section*{References}

\bibliographystyle{myjmva}
\bibliography{ref}

\begin{thebibliography}{30}
\expandafter\ifx\csname natexlab\endcsname\relax\def\natexlab#1{#1}\fi
\providecommand{\bibinfo}[2]{#2}
\ifx\xfnm\relax \def\xfnm[#1]{\unskip,\space#1}\fi
\bibitem[{Aitchison(1981)}]{Aitchison1981}
\bibinfo{author}{J.~Aitchison}, \bibinfo{title}{A new approach to null
  correlations of proportions}, \bibinfo{journal}{Journal of Mathematical
  Geology} \bibinfo{volume}{13} (\bibinfo{year}{1981})
  \bibinfo{pages}{175--189}.
\bibitem[{Bach(2008)}]{Bach2008}
\bibinfo{author}{F.~R. Bach}, \bibinfo{title}{Consistency of the group lasso
  and multiple kernel learning}, \bibinfo{journal}{Journal of Machine Learning
  Research} \bibinfo{volume}{9} (\bibinfo{year}{2008})
  \bibinfo{pages}{1179--1225}.
\bibitem[{Barab\'{a}si and Albert(1999)}]{Barabasi1999}
\bibinfo{author}{A.-L. Barab\'{a}si}, \bibinfo{author}{R.~Albert},
  \bibinfo{title}{Emergence of scaling in random networks},
  \bibinfo{journal}{Science} \bibinfo{volume}{286} (\bibinfo{year}{1999})
  \bibinfo{pages}{509--512}.
\bibitem[{Bickel and Levina(2008)}]{BicLev2008}
\bibinfo{author}{P.~J. Bickel}, \bibinfo{author}{E.~Levina},
  \bibinfo{title}{Regularized estimation of large covariance matrices},
  \bibinfo{journal}{Annals of Statistics} \bibinfo{volume}{36}
  (\bibinfo{year}{2008}) \bibinfo{pages}{199--227}.
\bibitem[{Charbonnier et~al.(2010)Charbonnier, Chiquet and
  Ambroise}]{Charbonnier2010}
\bibinfo{author}{C.~Charbonnier}, \bibinfo{author}{J.~Chiquet},
  \bibinfo{author}{C.~Ambroise}, \bibinfo{title}{Weighed-{L}asso for
  {S}tructured {N}etwork {I}nference from {T}ime {C}ourse {D}ata},
  \bibinfo{journal}{Statistical Applications in Genetics and Molecular Biology}
  \bibinfo{volume}{9} (\bibinfo{year}{2010}) \bibinfo{pages}{1544--6115}.
\bibitem[{Fan et~al.(2009)Fan, Feng and Wu}]{FanFengWu2009}
\bibinfo{author}{J.~Fan}, \bibinfo{author}{Y.~Feng}, \bibinfo{author}{Y.~Wu},
  \bibinfo{title}{Network exploration via the {A}daptive {L}asso and {SCAD}
  {P}enalties}, \bibinfo{journal}{The Annals of Applied Statistics}
  \bibinfo{volume}{3} (\bibinfo{year}{2009}) \bibinfo{pages}{521--541}.
\bibitem[{Fan and Li(2001)}]{FanLi2001}
\bibinfo{author}{J.~Fan}, \bibinfo{author}{R.~Li}, \bibinfo{title}{Variable
  selection via nonconcave penalized likelihood and its oracle properties.},
  \bibinfo{journal}{Journal of the American Statistical Association}
  \bibinfo{volume}{96} (\bibinfo{year}{2001}) \bibinfo{pages}{1348--1360}.
\bibitem[{Fan et~al.(2016)Fan, Liao and Liu}]{FanLiaoLiu2016}
\bibinfo{author}{J.~Fan}, \bibinfo{author}{Y.~Liao}, \bibinfo{author}{H.~Liu},
  \bibinfo{title}{An overview on the estimation of large covariance and
  precision matrices}, \bibinfo{journal}{Econometrics Journal}
  \bibinfo{volume}{19} (\bibinfo{year}{2016}) \bibinfo{pages}{1--32}.
\bibitem[{Faust et~al.(2015)Faust, Lahti, Gonze, de~Vos and Raes}]{Faust2015}
\bibinfo{author}{K.~Faust}, \bibinfo{author}{L.~Lahti},
  \bibinfo{author}{D.~Gonze}, \bibinfo{author}{W.~de~Vos},
  \bibinfo{author}{J.~Raes}, \bibinfo{title}{Metagenomics meets time series
  analysis: unraveling microbial community dynamics}, \bibinfo{journal}{Current
  Opinion in Microbiology} \bibinfo{volume}{25} (\bibinfo{year}{2015})
  \bibinfo{pages}{56--66}.
\bibitem[{Friedman and Alm(2012)}]{Friedman2012}
\bibinfo{author}{J.~Friedman}, \bibinfo{author}{E.~J. Alm},
  \bibinfo{title}{Inferring correlation networks from genomic survey data},
  \bibinfo{journal}{PLoS Comput Biol.} \bibinfo{volume}{8}
  (\bibinfo{year}{2012}) \bibinfo{pages}{e1002687}.
\bibitem[{Friedman et~al.(2008)Friedman, Hastie and
  Tibshirani}]{FriedmanHastieTibshirani2008}
\bibinfo{author}{J.~Friedman}, \bibinfo{author}{T.~Hastie},
  \bibinfo{author}{R.~Tibshirani}, \bibinfo{title}{Sparse inverse covariance
  estimation with the graphical lasso}, \bibinfo{journal}{Biostatistics}
  \bibinfo{volume}{9} (\bibinfo{year}{2008}) \bibinfo{pages}{432--441}.
\bibitem[{Gao et~al.(2012)Gao, Pu, Wu and Xu}]{Gao2012}
\bibinfo{author}{X.~Gao}, \bibinfo{author}{D.~Q. Pu}, \bibinfo{author}{Y.~Wu},
  \bibinfo{author}{H.~Xu}, \bibinfo{title}{Tuning {P}arameter {S}election for
  {P}enalized {L}ikelihood {E}stimation of {G}aussian {G}raphical {M}odel},
  \bibinfo{journal}{Statistica Sinica} \bibinfo{volume}{22}
  (\bibinfo{year}{2012}) \bibinfo{pages}{1123--1146}.
\bibitem[{Gilbert et~al.(2010)Gilbert, Meyer, Jansson, Gordon, Pace, Tiedje,
  Ley, Fierer, Field, Kyrpides et~al.}]{Gilbert2010}
\bibinfo{author}{J.~A. Gilbert}, \bibinfo{author}{F.~Meyer},
  \bibinfo{author}{J.~Jansson}, \bibinfo{author}{J.~Gordon},
  \bibinfo{author}{N.~Pace}, \bibinfo{author}{J.~Tiedje},
  \bibinfo{author}{R.~Ley}, \bibinfo{author}{N.~Fierer},
  \bibinfo{author}{D.~Field}, \bibinfo{author}{N.~Kyrpides}, et~al.,
  \bibinfo{title}{The earth microbiome project: Meeting {R}eport of the ``1st
  {EMP} meeting on {S}ample {S}election and {A}cquisition'' at {A}rgonne
  {N}ational {L}aboratory {O}ctober 6th 2010}, \bibinfo{journal}{Standards in
  {G}enomic {S}ciences} \bibinfo{volume}{3} (\bibinfo{year}{2010})
  \bibinfo{pages}{249}.
\bibitem[{van~der Heijden and Hartmann(2016)}]{HeiHart2016}
\bibinfo{author}{M.~G.~A. van~der Heijden}, \bibinfo{author}{M.~Hartmann},
  \bibinfo{title}{Networking in the plant microbiome}, \bibinfo{journal}{PLoS
  Biol.} \bibinfo{volume}{14} (\bibinfo{year}{2016}) \bibinfo{pages}{1--9}.
\bibitem[{Hero and Rajaratnam(2012)}]{HeroRaj2012}
\bibinfo{author}{A.~Hero}, \bibinfo{author}{B.~Rajaratnam}, \bibinfo{title}{Hub
  discovery in partial correlation graphs.}, \bibinfo{journal}{IEEE
  Transactions on Information Theory} \bibinfo{volume}{58}
  (\bibinfo{year}{2012}) \bibinfo{pages}{6064--6078}.
\bibitem[{Kurtz et~al.(2015)Kurtz, M\"{u}ller, Miraldi, Littman, Blaser and
  Bonneau}]{Kurtz2015}
\bibinfo{author}{Z.~D. Kurtz}, \bibinfo{author}{C.~L. M\"{u}ller},
  \bibinfo{author}{E.~R. Miraldi}, \bibinfo{author}{D.~R. Littman},
  \bibinfo{author}{M.~J. Blaser}, \bibinfo{author}{R.~A. Bonneau},
  \bibinfo{title}{Sparse and compositionally robust inference of microbial
  ecological networks}, \bibinfo{journal}{PLoS Comput Biol.}
  \bibinfo{volume}{11} (\bibinfo{year}{2015}) \bibinfo{pages}{e1004226}.
\bibitem[{Lam and Fan(2009)}]{LamFan2009}
\bibinfo{author}{C.~Lam}, \bibinfo{author}{J.~Fan},
  \bibinfo{title}{Sparsistency and rates of convergence in large covariance
  matrix estimation}, \bibinfo{journal}{Annals of Statistics}
  \bibinfo{volume}{37} (\bibinfo{year}{2009}) \bibinfo{pages}{4254--4278}.
\bibitem[{Li et~al.(2013)Li, Nasidze, Quinque, Li, Horz, Andr{\'e}, Garriga,
  Halbwax, Fischer and Stoneking}]{Li2013}
\bibinfo{author}{J.~Li}, \bibinfo{author}{I.~Nasidze},
  \bibinfo{author}{D.~Quinque}, \bibinfo{author}{M.~Li}, \bibinfo{author}{H.-P.
  Horz}, \bibinfo{author}{C.~Andr{\'e}}, \bibinfo{author}{R.~M. Garriga},
  \bibinfo{author}{M.~Halbwax}, \bibinfo{author}{A.~Fischer},
  \bibinfo{author}{M.~Stoneking}, \bibinfo{title}{The saliva microbiome of
  {P}an and {H}omo}, \bibinfo{journal}{BMC Microbiology} \bibinfo{volume}{13}
  (\bibinfo{year}{2013}) \bibinfo{pages}{204}.
\bibitem[{Liu and Ihler(2011)}]{LiuIhler2011}
\bibinfo{author}{Q.~Liu}, \bibinfo{author}{A.~Ihler}, \bibinfo{title}{Learning
  {S}cale {F}ree {N}etworks by {R}eweighted {L}1 {R}egularization},
  \bibinfo{journal}{Proceedings of the 14th {I}nternational {C}onference on
  {A}rtificial {I}ntelligence and {S}tatistics} \bibinfo{volume}{15}
  (\bibinfo{year}{2011}) \bibinfo{pages}{40--48}.
\bibitem[{Meinshausen and B\"{u}hlmann(2006)}]{MB2006}
\bibinfo{author}{N.~Meinshausen}, \bibinfo{author}{P.~B\"{u}hlmann},
  \bibinfo{title}{High dimensional graphs and variable selection with the
  lasso}, \bibinfo{journal}{Annals of Statistics} \bibinfo{volume}{34}
  (\bibinfo{year}{2006}) \bibinfo{pages}{1436--1462}.
\bibitem[{Mohan et~al.(2014)Mohan, London, Fazel, Witten and Lee}]{Mohan2014}
\bibinfo{author}{K.~Mohan}, \bibinfo{author}{P.~London},
  \bibinfo{author}{M.~Fazel}, \bibinfo{author}{D.~Witten},
  \bibinfo{author}{S.-I. Lee}, \bibinfo{title}{Node-based learning of multiple
  gaussian graphical models}, \bibinfo{journal}{The Journal of Machine Learning
  Research} \bibinfo{volume}{15} (\bibinfo{year}{2014})
  \bibinfo{pages}{445--488}.
\bibitem[{Negahban et~al.(2012)Negahban, Ravikumar, Wainwright and
  Yu}]{Negahban2012}
\bibinfo{author}{S.~N. Negahban}, \bibinfo{author}{P.~Ravikumar},
  \bibinfo{author}{M.~J. Wainwright}, \bibinfo{author}{B.~Yu},
  \bibinfo{title}{A unified framework for high-dimensional analysis of
  m-estimators with decomposable regularizers}, \bibinfo{journal}{Statist.
  Sci.} \bibinfo{volume}{27} (\bibinfo{year}{2012}) \bibinfo{pages}{538--557}.
\bibitem[{Obozinski et~al.(2011)Obozinski, Jacob and Vert}]{Obozinski2011}
\bibinfo{author}{G.~Obozinski}, \bibinfo{author}{L.~Jacob},
  \bibinfo{author}{J.-P. Vert}, \bibinfo{title}{Group lasso with overlaps: the
  latent group lasso approach}, \bibinfo{journal}{arXiv preprint
  arXiv:1110.0413}  (\bibinfo{year}{2011}).
\bibitem[{Rothman et~al.(2008)Rothman, Bickel, Levina and Zhu}]{Rothman2008}
\bibinfo{author}{A.~Rothman}, \bibinfo{author}{P.~J. Bickel},
  \bibinfo{author}{E.~Levina}, \bibinfo{author}{J.~Zhu}, \bibinfo{title}{Sparse
  permutation invariant covariance estimation}, \bibinfo{journal}{Electronic
  Journal of Statistics} \bibinfo{volume}{2} (\bibinfo{year}{2008})
  \bibinfo{pages}{494--515}.
\bibitem[{Shen et~al.(2012)Shen, Pan and Zhu}]{ShenPanZhu2012}
\bibinfo{author}{X.~Shen}, \bibinfo{author}{W.~Pan}, \bibinfo{author}{Y.~Zhu},
  \bibinfo{title}{Likelihood-based selection and sharp parameter estimation},
  \bibinfo{journal}{Journal of the American Statistical Association}
  \bibinfo{volume}{107} (\bibinfo{year}{2012}) \bibinfo{pages}{223--232}.
\bibitem[{Tan et~al.(2014)Tan, London, Mohan, Lee, Fazel and Witten}]{Tan2014}
\bibinfo{author}{K.~M. Tan}, \bibinfo{author}{P.~London},
  \bibinfo{author}{K.~Mohan}, \bibinfo{author}{S.-I. Lee},
  \bibinfo{author}{M.~Fazel}, \bibinfo{author}{D.~Witten},
  \bibinfo{title}{Learning {G}raphical {M}odels with {H}ubs},
  \bibinfo{journal}{Journal of Machine Learning Research} \bibinfo{volume}{15}
  (\bibinfo{year}{2014}) \bibinfo{pages}{3297--3331}.
\bibitem[{Turnbaugh et~al.(2007)Turnbaugh, Ley, Hamady, Fraser-Liggett, Knight
  and Gordon}]{Turnbaugh2007}
\bibinfo{author}{P.~J. Turnbaugh}, \bibinfo{author}{R.~E. Ley},
  \bibinfo{author}{M.~Hamady}, \bibinfo{author}{C.~M. Fraser-Liggett},
  \bibinfo{author}{R.~Knight}, \bibinfo{author}{J.~I. Gordon},
  \bibinfo{title}{The {H}uman {M}icrobiome {P}roject},
  \bibinfo{journal}{Nature} \bibinfo{volume}{449} (\bibinfo{year}{2007})
  \bibinfo{pages}{804--810}.
\bibitem[{Witten et~al.(2011)Witten, Friedman and Simon}]{Witten2011}
\bibinfo{author}{D.~M. Witten}, \bibinfo{author}{J.~H. Friedman},
  \bibinfo{author}{N.~Simon}, \bibinfo{title}{New insights and faster
  computations for the graphical lasso}, \bibinfo{journal}{Journal of
  Computational and Graphical Statistics} \bibinfo{volume}{20}
  (\bibinfo{year}{2011}) \bibinfo{pages}{892--900}.
\bibitem[{Yuan and Lin(2007)}]{YuanLin2007}
\bibinfo{author}{M.~Yuan}, \bibinfo{author}{Y.~Lin}, \bibinfo{title}{Model
  selection and estimation in the gaussian graphical model},
  \bibinfo{journal}{Biometrika} \bibinfo{volume}{94} (\bibinfo{year}{2007})
  \bibinfo{pages}{19--35}.
\bibitem[{Zou(2006)}]{Zou2006}
\bibinfo{author}{H.~Zou}, \bibinfo{title}{The adaptive lasso and its oracle
  properties}, \bibinfo{journal}{Journal of the American Statistical
  Association} \bibinfo{volume}{101} (\bibinfo{year}{2006})
  \bibinfo{pages}{1418--1429}.

\end{thebibliography}

\renewcommand{\baselinestretch}{1}

\setlength{\tabcolsep}{1.8pt}
\begin{table}[p]
\centering
\footnotesize
\begin{tabular}{ccccccccccc}
  \hline
  \hline
Method  &  True Pos. & True Neg.  & Perc. of Correctly &  Perc. of Correctly & Number of & Frobenius  \\ 
& Rate &Rate & Estimated Hub & Estimated Hub / & Estimated & Norm \\ 
 &  (TPR) & (TNR)  &   Edges &   Non-Hub Nodes  & Edges  \\
 \hline 
   \hline
   \emph{Simulation (i)} &         &   \\
       &         & \multicolumn{3}{c}{$n=100, p=50$} \\
\verb+glasso+ & 72.69 (0.26) & 84.03 (0.51) & 61.27 (0.41) & 100 (0)/32.85 (2.01) & 234.71 (6.01) & 3.30 (0.02) \\ 
\verb+Ada-glasso+ &  80.51 (0.55) & 96.40 (0.24) & 74.51 (0.86) & 100 (0)/88.65 (1.38) & 105.87 (3.37) & 1.73 (0.02) \\ 
\verb+SF+ &  72.98 (0.27) & 95.56 (0.18) & 62.95 (0.42) & 100 (0)/75.48 (1.09) & 104.57 (2.25) & 2.22 (0.02) \\ 
\verb+HGL+ ($c=0.50$) &  74.05 (0.22) & 82.74 (0.40) & 63.53 (0.36) & 100 (0)/32.60 (1.53) & 251.26 (4.73) & 3.24 (0.02) \\ 
\verb+HGL+ ($c=0.75$) &  73.60 (0.19) & 84.20 (0.24) & 62.98 (0.32) & 100 (0)/39.31 (0.80) & 234.06 (2.82) & 3.31 (0.02) \\ 
\verb+hw.glasso+ &  87.06 (0.37) & 98.67 (0.09) & 85.85 (0.61) & 100 (0)/99.33 (0.16) & 89.52 (1.34) & 1.14 (0.02) \\ 
2-step \verb+hw.glasso+ &  94.27 (0.16) & 98.49 (0.15)& 98.51 (0.28) & 100 (0)/99.33 (0.16) & 101.87 (1.59) & 0.94 (0.03) \\ 
2-step \verb+hw.glasso+  &   94.57 (0.08) & 99.20 (0.02) & 99.08 (0.13) & 100 (0)/100 (0) & 94.29 (0.26) & 0.79 (0.01) \\ 
(known hubs) & \\
\hline
  &  & \multicolumn{3}{c}{$n=100, p=100$} \\
\verb+glasso+ &  48.10 (0.26) & 94.40 (0.28) & 38.28 (0.33) & 99.50 (0.35) / 73.70 (1.60) & 384.73 (13.73) & 7.29 (0.05) \\ 
\verb+Ada-glasso+ &  58.31 (0.19) & 96.55 (0.03) & 52.97 (0.26) & 100 (0) / 99.24 (0.10) & 334.78 (1.42) & 4.49 (0.02) \\  
\verb+SF+ & 53.08 (0.33) & 97.94 (0.07) & 46.53 (0.46) & 99.25 (0.56) / 95.05 (0.37) & 246.59 (4.50) & 5.34 (0.04) \\ 
\verb+HGL+ ($c=0.50$) &  56.12 (0.17) & 84.26 (0.29) & 47.45 (0.20) & 100 (0) / 19.91 (1.51) & 886.60 (13.72) & 6.43 (0.02) \\
\verb+HGL+ ($c=0.75$) &  50.81 (0.31) & 92.82 (0.33) & 42.05 (0.40) & 99.50 (0.50) / 65.26 (1.67) & 469.76 (16.53) & 7.30 (0.04) \\   
\verb+hw.glasso+  &   70.55 (0.49) & 99.60 (0.01) & 72.77 (0.72) & 100 (0) / 100 (0) & 253.23 (2.68) & 2.75 (0.03) \\ 
2-step \verb+hw.glasso+  &  79.24 (0.36) & 99.23 (0.01) & 85.56 (0.52) & 100 (0) / 100 (0) & 311.58 (2.17) & 2.62 (0.03) \\
2-step \verb+hw.glasso+   &  79.24 (0.36) & 99.23 (0.01) & 85.56 (0.52) & 100 (0) / 100 (0) & 311.58 (2.17) & 2.62 (0.03) \\
(known hubs) & \\
  \hline
&  & \multicolumn{3}{c}{$n=100, p=200$} \\
\verb+glasso+ & 24.76 (0.22) & 99.30 (0.03) & 16.01 (0.28) & 66.38 (1.11) / 99.18 (0.11) & 336.06 (9.93) & 14.98 (0.09) \\ 
\verb+Ada-glasso+ &  27.30 (0.13) & 99.00 (0.03) & 19.25 (0.16) & 78.75 (0.91) / 99.99 (0.01) & 432.65 (6.94) & 13.28 (0.06) \\ 
\verb+SF+ &  28.54 (0.14) & 99.50 (0.02) & 20.98 (0.17) & 68.12 (0.86) / 99.81 (0.03) & 361.03 (5.05) & 11.19 (0.05) \\ 
\verb+HGL+ ($c=0.50$) &  49.69 (0.24) & 58.26 (0.35) & 41.73 (0.26) & 100 (0) / 0 (0) & 8319.62 (68.81) & 73.08 (0.97) \\
\verb+HGL+ ($c=0.75$) &  33.69 (0.09) & 92.95 (0.21) & 26.28 (0.10) & 93.12 (0.62) / 69.73 (1.29) & 1654.97 (38.89) & 13.14 (0.05) \\  
\verb+hw.glasso+  & 31.18 (0.16) & 99.81 (0.01) & 24.53 (0.21) & 83.62 (1.11) / 100 (0) & 347.32 (3.47) & 8.91 (0.03) \\
2-step \verb+hw.glasso+  & 42.15 (0.28) & 99.69 (0.001) & 38.75 (0.36) & 83.62 (1.11) / 100 (0) & 548.95 (5.10) & 8.76 (0.04) \\
2-step \verb+hw.glasso+  & 45.18 (0.18) & 99.65 (0.005) & 42.67 (0.23) & 100 (0) / 100 (0)  & 605.59 (3.38) & 8.76 (0.04) \\ 
(known hubs) & \\
 \hline
&  & \multicolumn{3}{c}{$n=250, p=500$} \\
\verb+glasso+ & 14.90(0.07) & 99.44(0.01) & 11.51(0.08) & 47.45(0.45)/99.90(0.02) & 1565.71(19.16) & 32.50(0.07) \\ 
\verb+Ada-glasso+ & 16.47(0.03) & 99.45(0.003) & 13.34(0.03) & 60.20(0.51)/100(0) & 1703.01(4.43) & 27.27(0.04) \\ 
\verb+SF+ & 18.27(0.07) & 99.72(0.004) & 15.65(0.08) & 65.20(0.43)/100(0) & 1567.36(10.69) & 25.71(0.04) \\ 
\verb+HGL+ ($c=0.50$) & 31.00(0.05) & 82.85(0.02) & 27.34(0.05) & 100(0)/0(0) & 22294.14(26.20) & 45.76(0.05) \\ 
\verb+HGL+ ($c=0.75$) & 22.41(0.27) & 92.56(0.22) & 18.77(0.27) & 98.10(0.43)/ 80.35(1.29) & 10249.74(281.84) & 31.64(0.24) \\ 
\verb+hw.glasso+  & 21.02(0.12) & 99.90(0.002) & 19.03(0.15) & 77.95(0.88)/100(0) & 1612.38(13.02) & 21.26(0.06) \\ 
2-step \verb+hw.glasso+  & 26.16(0.15) & 99.81(0.002) & 25.26(0.18) & 78.05(0.88)/100(0) & 2209.02(15.08) & 23.09(0.06) \\ 
2-step \verb+hw.glasso+  & 27.70(0.12) & 99.80(0.002) & 27.12(0.14) & 100(0)/100(0) & 2369.31(12.50) & 23.60(0.07) \\ 
(known hubs) & \\
 \hline
&  & \multicolumn{3}{c}{$n=1000, p=500$} \\
\verb+glasso+ & 26.64(0.05) & 97.98(0.02) & 24.90(0.06) & 98.90(0.21)/97.44(0.04) & 4380.73(24.70) & 25.67(0.03) \\  
\verb+Ada-glasso+ & 28.48(0.27) & 99.36(0.03) & 27.48(0.30) & 100(0)/100(0) & 2955.13(59.37) & 20.21(0.07) \\ 
\verb+SF+ & 37.59(0.04) & 99.32(0.004) & 38.65(0.05) & 100(0)/99.99(0.01) & 3870.29(6.83) & 18.71(0.01) \\ 
\verb+HGL+ ($c=0.50$) & 35.25(0.04) & 93.68(0.01)& 33.30(0.04) & 100(0) /87.89(0.12) & 10171.54(11.76) & 22.98(0.02) \\ 
\verb+HGL+ ($c=0.75$) & 31.06(0.09) & 96.05(0.05) & 29.29(0.08) & 100(0)/95.11(0.11) & 7033.16(68.45) & 23.81(0.03) \\ 
\verb+hw.glasso+  & 43.02(0.13) & 99.68(0.003) & 45.46(0.15) & 100(0)/100(0) & 3968.17(14.70) & 15.00(0.03) \\ 
2-step \verb+hw.glasso+  & 53.30(0.23) & 99.69(0.002) & 58.03(0.28) & 100(0)/100(0) & 4929.58(24.09) & 15.04(0.07) \\ 
2-step \verb+hw.glasso+  & 53.30(0.23) & 99.69(0.002) & 58.03(0.28) & 100(0)/100(0) & 4929.58(24.09) & 15.04(0.07) \\  
(known hubs) & \\
   \hline
\end{tabular}
\caption{Means (and standard errors) of different performance measures over 100 replications for the graphical lasso (glasso), adaptive graphical lasso (Ada-glasso), scale-free network approach (SF), hubs graphical lasso (HGL), hubs weighted graphical lasso (hw.glasso), and 2-step hw.glasso. } 
\label{tab:Simulation1}
\end{table}

\renewcommand{\baselinestretch}{1.0}

\setlength{\tabcolsep}{1.8pt}
\begin{table}[th] 
\footnotesize
\centering
\begin{tabular}{ccccccccccc}
  \hline
  \hline
Method  & True Pos. & True Neg.  & Perc. of Correctly &  Perc. of Correctly & Number of & Frobenius  \\ 
&  Rate &Rate & Estimated Hub & Estimated Hub / & Estimated & Norm \\ 
 &   (TPR) & (TNR)   &   Edges &   Non-Hub Nodes  & Edges  \\
 \hline 
   \hline
   \emph{Simulation (ii)} &         &  \\
  &         & \multicolumn{3}{c}{$n=100, p=50$} \\
\verb+glasso+ &   90.42 (0.25) & 93.48 (0.15) & 88.64 (0.31) & 100 (0)/69.60 (0.65) & 107.74 (1.87) & 1.01 (0.01) \\ 
\verb+Ada-glasso+ & 91.21 (0.27) & 98.15 (0.11) & 93.00 (0.50) & 100 (0)/99.17 (0.14) & 53.18 (1.45) & 0.51 (0.01) \\ 
\verb+SF+ &   87.83 (0.22) & 97.59 (0.08) & 90.43 (0.32) & 100 (0)/92.31 (0.52) & 56.79 (1.05) & 0.65 (0.01) \\ 
\verb+HGL+ ($c=0.50$) & 91.27 (0.24) & 92.17 (0.24) & 89.75 (0.31) & 100 (0)/64.83 (1.05) & 124.10 (2.96) & 1.01 (0.01) \\ 
\verb+HGL+ ($c=0.75$) &  90.34 (0.23) & 93.19 (0.13) & 89.29 (0.29) & 100 (0)/68.77 (0.63) & 111.15 (1.57) & 1.04 (0.01) \\ 
\verb+hw.glasso+ & 91.48 (0.25) & 98.47 (0.07) & 94.68 (0.40) & 100 (0)/99.50 (0.10) & 49.55 (0.96) & 0.46 (0.01) \\ 
2-step \verb+hw.glasso+ &  87.65 (0.17) & 96.92 (0.07) & 96.68 (0.31) & 100 (0)/99.50 (0.10) & 64.59 (0.81) & 0.52 (0.01) \\ 
2-step \verb+hw.glasso+ & 87.21 (0.15) & 97.17 (0.05) & 96.57 (0.28) & 100 (0)/100 (0) & 61.23 (0.65) & 0.51 (0.01) \\ 
 (known hubs) & \\
\hline
&  & \multicolumn{3}{c}{$n=100, p=100$} \\
\verb+glasso+ &   66.17 (0.36) & 97.48 (0.07) & 57.84 (0.62) & 99.25 (0.43) / 92.75 (0.35) & 198.98 (4.25) & 2.54 (0.01) \\ 
\verb+Ada-glasso+ & 72.86 (0.19) & 98.57 (0.02) & 72.56 (0.37) & 99.75 (0.25) / 100 (0) & 164.85 (0.92) & 1.59 (0.01) \\  
\verb+SF+ &   71.31 (0.26) & 98.38 (0.04) & 72.86 (0.47) & 100 (0) / 97.74 (0.15) & 169.79 (2.34) & 1.77 (0.01) \\ 
\verb+HGL+ ($c=0.50$) &  74.01 (0.27) & 94.27 (0.16) & 69.57 (0.39) & 100 (0) / 76.81 (0.84) & 373.57 (8.07) & 2.32 (0.01) \\ 
\verb+HGL+ ($c=0.75$) &  68.34 (0.33) & 96.87 (0.09) & 62.33 (0.54) & 100 (0) / 88.56 (0.49) & 234.09 (5.24) & 2.52 (0.01) \\  
\verb+hw.glasso+  &  74.94 (0.24) & 99.11 (0.03) & 83.49 (0.42) & 100 (0) / 100 (0) & 144.86 (1.77) & 1.25 (0.01) \\ 
2-step \verb+hw.glasso+  &  75.02 (0.16) & 97.83 (0.03) & 88.98 (0.36) & 100 (0) / 100 (0) & 206.12 (1.93) & 1.44 (0.01) \\
2-step \verb+hw.glasso+ & 75.02 (0.16) & 97.83 (0.03) & 88.98 (0.36) & 100 (0) / 100 (0) & 206.12 (1.93) & 1.44 (0.01) \\ 
 (known hubs)  & \\
\hline
&  & \multicolumn{3}{c}{$n=100, p=200$} \\
\verb+glasso+ &  36.77 (0.24) & 99.47 (0.02) & 23.10 (0.38) & 48.00 (1.27) / 99.79 (0.03) & 222.39 (6.03) & 5.71 (0.02) \\ 
\verb+Ada-glasso+ &  41.83 (0.21) & 99.30 (0.03) & 31.27 (0.33) & 60.25 (1.26) / 100 (0) & 299.40 (6.88) & 5.17 (0.02) \\ 
\verb+SF+ &   43.25 (0.22) & 99.39 (0.02) & 34.32 (0.35) & 68.38 (1.01) / 99.72 (0.03) & 294.44 (4.61) & 4.47 (0.02) \\ 
\verb+HGL+ ($c=0.50$) &  73.53 (0.21) & 57.97 (0.31) & 66.61 (0.27) & 100 (0) / 0 (0) & 8522.47 (61.84) & 30.42 (0.41) \\ 
\verb+HGL+ ($c=0.75$) & 49.11 (0.22) & 96.56 (0.10) & 41.35 (0.28) & 92.25 (0.94) / 92.83 (0.43) & 888.32 (21.62) & 5.13 (0.02) \\ 
\verb+hw.glasso+ & 50.98 (0.29) & 99.51 (0.01) & 47.95 (0.47) & 85.88 (0.98) / 100 (0) & 338.71 (4.53) & 3.43 (0.02)\\
2-step \verb+hw.glasso+  &  56.55 (0.30) & 98.95 (0.02) & 58.04 (0.53) & 85.88 (0.98) / 100 (0) & 493.88 (5.16) & 3.75 (0.02) \\
2-step \verb+hw.glasso+  &  58.80 (0.26) & 98.92 (0.02) & 61.92 (0.46) & 100 (0) / 100 (0) & 519.27 (5.20) & 3.76 (0.02) \\ 
(known hubs) & \\
   \hline
&  & \multicolumn{3}{c}{$n=250, p=500$} \\
\verb+glasso+ & 20.90(0.10) & 99.67(0.01) & 14.95(0.15) & 26.60(0.52)/100(0) & 868.82(14.55) & 12.21(0.02) \\ 
\verb+Ada-glasso+ &  24.21(0.05) & 99.67(0.003) & 19.67(0.07) & 26.95(0.50)/100(0) & 1019.44(4.47) & 10.70(0.01) \\  
\verb+SF+ & 23.53(0.11) & 99.73(0.00) & 19.13(0.16) & 37.35(0.41)/100(0) & 918.31(9.45) & 10.42(0.02) \\ 
\verb+HGL+ ($c=0.50$) &  51.18(0.14) & 81.82(0.11) & 47.59(0.16) & 100(0)/0(0) & 23805.59(143.13) & 20.07(0.11) \\ 
\verb+HGL+ ($c=0.75$) & 30.10(0.23) & 96.68(0.08) & 24.62(0.27) & 47.85(2.67)/99.99(0.01) & 4896.21(111.53) & 11.08(0.03) \\ 
\verb+hw.glasso+  & 29.03(0.17) & 99.77(0.004) & 27.53(0.25) & 44.35(0.84)/100(0) & 1128.74(11.86) & 8.82(0.02) \\
2-step \verb+hw.glasso+  & 28.14(0.21) & 99.53(0.00) & 26.29(0.32) & 44.55(0.84)/100(0) & 1375.00(14.88) & 9.67(0.02) \\ 
2-step \verb+hw.glasso+  & 34.14(0.20) & 99.51(0.01) & 35.36(0.30) & 90.15(1.11)/ 100(0) & 1679.77(16.34) & 9.93(0.03) \\ 
(known hubs) & \\
   \hline
&  & \multicolumn{3}{c}{$n=1000, p=500$} \\
\verb+glasso+ & 41.41(0.10) & 98.89(0.01) & 42.21(0.12) & 91.35(0.43)/99.86(0.01) & 2760.33(14.86) & 9.48(0.01) \\ 
\verb+Ada-glasso+ & 47.54(0.53) & 99.17(0.03) & 50.53(0.63) & 97.20(0.38)/100(0) & 2710.44(66.12) & 7.17(0.05) \\ 
\verb+SF+ & 54.39(0.06) & 99.29(0.004) & 63.60(0.09) & 99.80(0.10)/100(0) & 2882.43(5.87) & 6.96(0.01) \\ 
\verb+HGL+ ($c=0.50$) &  44.71(0.34) & 98.48(0.05) & 45.82(0.36) & 95.75(0.33)/99.61(0.03) & 3411.80(72.09) & 9.18(0.03) \\ 
\verb+HGL+ ($c=0.75$) & 31.82(0.40) & 99.51(0.02) & 30.46(0.52) & 53.45(1.46)/99.99(0.01) & 1568.98(43.07) & 10.62(0.05) \\ 
\verb+hw.glasso+  & 56.37(0.14) & 99.30(0.01) & 66.86(0.18) & 99.80(0.10)/100(0) & 2963.33(13.17) & 5.63(0.02) \\ 
2-step \verb+hw.glasso+  & 67.27(0.16) & 98.91(0.01) & 85.47(0.24) & 99.80(0.10)/100(0) & 3933.57(19.10) & 5.39(0.03) \\ 
2-step \verb+hw.glasso+  & 67.32(0.16) & 98.91(0.01) & 85.55(0.24) & 100(0)/100(0) & 3935.49(19.36) & 5.40(0.03) \\ 
(known hubs) & \\
    \hline
\end{tabular}
\caption{Means (and standard errors) of different performance measures over 100 replications for the graphical lasso (glasso), adaptive graphical lasso (Ada-glasso), scale-free network approach (SF), hubs graphical lasso (HGL), hubs weighted graphical lasso (hw.glasso), and 2-step hw.glasso. } 
\label{tab:Simulation2}
\end{table}

\renewcommand{\baselinestretch}{1.0}

\setlength{\tabcolsep}{2.5pt}
\begin{table}[th] 
\footnotesize
\centering
\begin{tabular}{cccccccccc}
  \hline
     \hline
Method  & True Pos. & True Neg.  & Perc. of Correctly &  Perc. of Correctly & Number of & Frobenius  \\ 
&  Rate &Rate & Estimated Hub & Estimated Hub / & Estimated & Norm \\ 
 &   (TPR) & (TNR)   &   Edges &   Non-Hub Nodes  & Edges  \\
 \hline 
   \hline
   \emph{Simulation (iii)} &         &   \\
 &         & \multicolumn{3}{c}{$n=100, p=50$} \\
\verb+glasso+ &  90.43 (0.27) & 89.16 (0.25) & 95.03 (0.39) & 100 (0)/46.44 (1.33) & 178.92 (3.12) & 1.52 (0.01) \\ 
\verb+Ada-glasso+ &  87.93 (0.23) & 97.21 (0.04) & 95.74 (0.32) & 100 (0)/96.58 (0.28) & 82.67 (0.58) & 0.80 (0.01) \\ 
\verb+SF+ & 86.44 (0.27) & 95.93 (0.11) & 95.47 (0.37) & 100 (0)/79.56 (0.80) & 95.76 (1.54) & 1.03 (0.01) \\ 
\verb+HGL+ ($c=0.50$) &   89.61 (0.25) & 89.39 (0.24) & 94.74 (0.37) & 100 (0)/50.40 (1.13) & 175.31 (2.94) & 1.60 (0.01) \\
\verb+HGL+ ($c=0.75$) &  88.78 (0.20) & 90.47 (0.13) & 94.37 (0.36) & 100 (0)/56.21 (0.56) & 161.84 (1.66) & 1.67 (0.01) \\ 
\verb+hw.glasso+  &  87.30 (0.31) & 97.67 (0.11) & 96.21 (0.32) & 100 (0)/96.92 (0.45) & 76.53 (1.53) & 0.78 (0.01) \\ 
2-step \verb+hw.glasso+  &  79.92 (0.33) & 95.61 (0.20) & 98.95 (0.18) & 100 (0)/96.15 (0.64) & 92.09 (2.63) & 0.94 (0.01) \\ 
2-step \verb+hw.glasso+  & 78.03 (0.16) & 97.14 (0.05) & 99.92 (0.05) & 100 (0)/100 (0) & 72.16 (0.65) & 0.88 (0.01) \\ 
(known hubs) & \\
\hline
 &         & \multicolumn{3}{c}{$n=100, p=100$} \\
\verb+glasso+ &  51.91 (0.30) & 97.41 (0.06) & 42.35 (0.46) & 66.50 (1.95) / 87.89 (0.37) & 199.94 (3.78) & 3.57 (0.02) \\ 
\verb+Ada-glasso+ & 49.45 (0.33) & 99.59 (0.03) & 43.01 (0.58) & 68.50 (1.80) / 99.99 (0.01) & 88.70 (2.18) & 2.54 (0.01) \\ 
\verb+SF+ &  49.10 (0.41) & 98.76 (0.04) & 39.98 (0.80) & 64.50 (1.92) / 97.00 (0.22) & 126.91 (3.14) & 2.79 (0.02) \\ 
\verb+HGL+ ($c=0.50$) &  57.68 (0.29) & 95.55 (0.10) & 53.82 (0.46) & 95.25 (1.11) / 79.58 (0.48) & 307.45 (5.35) & 3.45 (0.01) \\
\verb+HGL+ ($c=0.75$) & 52.01 (0.49) & 96.90 (0.10) & 43.84 (0.82) & 75.75 (2.29) / 85.84 (0.47) & 224.18 (6.07) & 3.64 (0.02) \\ 
\verb+hw.glasso+  &  57.97 (0.34) & 99.31 (0.03) & 63.07 (0.63) & 95.25 (0.99) / 99.99 (0.01) & 131.19 (2.34) & 2.10 (0.01) \\ 
2-step \verb+hw.glasso+  &   61.46 (0.35) & 98.77 (0.04) & 75.27 (0.82) & 95.25 (0.99) / 99.99 (0.01) & 168.87 (2.68) & 2.31 (0.02) \\
2-step \verb+hw.glasso+  &  62.65 (0.26) & 98.76 (0.04) & 78.07 (0.61) & 100 (0) / 100 (0) & 173.14 (2.53) & 2.30 (0.02) \\ 
(known hubs) & \\
\hline
 &         & \multicolumn{3}{c}{$n=100, p=200$} \\
\verb+glasso+ & 33.57 (0.23) & 99.04 (0.04) & 30.07 (0.36) & 65.38 (0.75) / 99.74 (0.05) & 392.68 (9.40) & 8.87 (0.03) \\ 
\verb+Ada-glasso+ &  36.13 (0.14) & 99.13 (0.03) & 34.70 (0.21) & 69.12 (1.07) / 100 (0) & 405.32 (6.17) & 6.99 (0.02) \\ 
\verb+SF+ &  35.82 (0.14) & 99.46 (0.02) & 35.89 (0.23) & 67.25 (0.68) / 99.97 (0.01) & 340.27 (4.27) & 6.73 (0.02) \\ 
\verb+HGL+ ($c=0.50$) &   42.27 (0.10) & 95.90 (0.09) & 40.85 (0.20) & 82.00 (0.74) / 92.01 (0.32) & 1091.75 (17.25) & 7.87 (0.02) \\ 
\verb+HGL+ ($c=0.75$) & 42.34 (0.07) & 95.69 (0.03) & 40.50 (0.12) & 82.25 (0.69) / 91.79 (0.29) & 1130.89 (6.53) & 7.82 (0.02) \\ 
\verb+hw.glasso+  &   37.78 (0.17) & 99.66 (0.01) & 41.03 (0.32) & 76.62 (1.25) / 100 (0) & 325.78 (3.59) & 5.37 (0.02) \\ 
2-step \verb+hw.glasso+&   43.17 (0.30) & 99.26 (0.02) & 51.72 (0.58) & 76.62 (1.25) / 100 (0) & 467.41 (5.80) & 5.72 (0.03) \\ 
2-step \verb+hw.glasso+ &   47.59 (0.15) & 99.25 (0.01) & 60.26 (0.29) & 100 (0) / 100 (0) & 524.46 (4.29) & 5.68 (0.03) \\ 
(known hubs) & \\
 \hline
  &         & \multicolumn{3}{c}{$n=250, p=500$} \\
\verb+glasso+ &  16.00(0.09) & 99.52(0.01) & 14.40(0.14) & 31.10(0.39)/99.98(0.01) & 1186.32(18.43) & 18.96(0.03) \\ 
\verb+Ada-glasso+ &  18.10(0.09) & 99.69(0.01) & 18.32(0.13) & 28.90(0.74)/100(0) & 1131.18(13.40) & 15.86(0.02) \\ 
\verb+SF+ &  18.14(0.16) & 99.74(0.01) & 18.84(0.27) & 43.95(0.76)/ 100(0) & 1084.95(17.17) & 16.20(0.03) \\ 
\verb+HGL+ ($c=0.50$) &  36.74(0.05) & 90.44(0.01) & 36.02(0.07) & 100(0)/40.79(0.40) & 13383.04(14.94) & 18.75(0.01) \\ 
\verb+HGL+ ($c=0.75$) & 

 30.26(0.08) & 95.58(0.05) & 30.85(0.08) & 96.80(0.39)/99.86(0.02) & 6846.39(61.58) & 16.90(0.02) \\ 
 
\verb+hw.glasso+  & 23.18(0.13) & 99.83(0.003) & 28.29(0.22) & 66.70(1.12)/100(0) & 1328.33(11.95) & 13.75(0.03) \\ 
2-step \verb+hw.glasso+&  24.75(0.20) & 99.74(0.004) & 31.23(0.36) & 66.85(1.13)/100(0) & 1548.12(15.75) & 15.17(0.03) \\ 
2-step \verb+hw.glasso+ &  
27.93(0.15) & 99.77(0.01) & 36.87(0.27) & 97.60(0.34)/100(0) & 1731.88(15.64) & 15.34(0.04) \\ 
(known hubs) & \\
  \hline
  &         & \multicolumn{3}{c}{$n=1000, p=500$} \\
  \verb+glasso+ & 36.04(0.09) & 98.39(0.01) & 41.07(0.09) & 97.00(0.25)/99.09(0.03) & 3928.19(18.43) & 14.83(0.01) \\ 
\verb+Ada-glasso+ &   42.19(0.24) & 98.86(0.02) & 51.24(0.28) & 100(0)/100(0) & 3803.20(38.99) & 11.22(0.04) \\ 
\verb+SF+ &   43.78(0.24) & 99.29(0.01) & 59.94(0.32) & 100(0)/100(0) & 3404.85(33.92) & 11.58(0.04) \\ 
\verb+HGL+ ($c=0.50$) & 41.75(0.04) & 97.62(0.01) & 46.76(0.06) & 100(0)/ 98.70(0.04) & 5246.92(8.32) & 13.93(0.01) \\ 
\verb+HGL+ ($c=0.75$) & 

 32.73(0.26) & 98.71(0.02) & 37.46(0.29) & 91.00(0.71)/99.22(0.03) & 3322.02(46.36) & 15.37(0.04) \\ 

\verb+hw.glasso+  &   44.80(0.12) & 99.32(0.01) & 63.25(0.16) & 100(0)/100(0) & 3440.72(16.25) & 9.80(0.02) \\ 
2-step \verb+hw.glasso+&    52.58(0.11) & 99.19(0.01) & 80.55(0.20) & 100(0)/100(0) & 4150.51(17.23) & 9.47(0.03) \\ 
2-step \verb+hw.glasso+ &  52.58(0.11) & 99.19(0.01)& 80.55(0.20) & 100(0)/100(0) & 4150.51(17.23) & 9.47(0.03) \\ 
(known hubs) & \\
  \hline
 \\
\end{tabular}
\caption{Means (and standard errors) of different performance measures over 100 replications for the graphical lasso (glasso), adaptive graphical lasso (Ada-glasso), scale-free network approach (SF), hubs graphical lasso (HGL), hubs weighted graphical lasso (hw.glasso), and 2-step hw.glasso.
} 
\label{tab:Simulation3}
\end{table}

\renewcommand{\baselinestretch}{1.0}

\setlength{\tabcolsep}{2.5pt}
\begin{table}[th] 
\footnotesize
\centering
\begin{tabular}{cccccccccc}
  \hline
    \hline
Method  & True Pos. & True Neg.  & Perc. of Correctly &  Perc. of Correctly & Number of & Frobenius  \\ 
& Rate &Rate & Estimated Hub & Estimated Hub / & Estimated & Norm \\ 
 &   (TPR) & (TNR)   &   Edges &   Non-Hub Nodes  & Edges  \\
 \hline 
   \hline
   \emph{Simulation (iv)} &         & \\
    &         & \multicolumn{3}{c}{$n=100, p=50$} \\
\verb+glasso+ &  88.83 (0.41) & 90.34 (0.30) & 95.62 (0.42) & 100 (0)/56.75 (1.26) & 151.50 (3.83) & 1.28 (0.01) \\ 
\verb+Ada-glasso+ &  89.79 (0.29) & 97.51 (0.08) & 95.31 (0.37) & 100 (0)/98.85 (0.16) & 68.15 (1.05) & 0.60 (0.01) \\ 
\verb+SF+ &  83.71 (0.27) & 97.10 (0.10) & 95.31 (0.47) & 99.50 (0.50)/93.79 (0.60) & 67.01 (1.33) & 0.83 (0.01) \\
\verb+hw.glasso+  &   87.22 (0.30) & 98.52 (0.07) & 95.22 (0.41) & 100 (0)/99.73 (0.08) & 53.74 (0.99) & 0.55 (0.01) \\ 
2-step \verb+hw.glasso+  &  83.84 (0.18) & 96.51 (0.06) & 99.31 (0.13) & 100 (0)/99.71 (0.08) & 74.10 (0.80) & 0.62 (0.01) \\ 
2-step \verb+hw.glasso+  &  83.32 (0.16) & 96.64 (0.05) & 99.12 (0.14) & 100 (0)/100 (0) & 71.97 (0.63) & 0.62 (0.01) \\ 
(known hubs) & \\
\hline
 &         & \multicolumn{3}{c}{$n=100, p=100$} \\
\verb+glasso+ & 70.58 (0.31) & 98.34 (0.08) &  64.92 (0.84) & 83.00 (2.39) / 97.04 (0.43) & 121.01 (4.38) & 1.57 (0.01) \\  
\verb+Ada-glasso+ &   78.69 (0.35) & 98.85 (0.06) & 80.29 (0.71) & 79.00 (2.71) / 99.99 (0.01) & 112.62 (3.28) & 0.91 (0.01) \\  
\verb+SF+ & 72.25 (0.32) & 99.27 (0.03) &  72.19 (0.98) & 71.00 (2.67) / 99.96 (0.03) & 79.14 (2.04) & 1.10 (0.01) \\ 
\verb+hw.glasso+ & 76.94 (0.27) & 99.29 (0.03) & 80.53 (0.67) & 80.00 (2.22) / 100 (0) & 87.66 (1.64) & 1.03 (0.01) \\ 
2-step \verb+hw.glasso+  &  75.01 (0.46) & 98.29 (0.03) & 83.10 (1.55) & 80.33 (2.23) / 100 (0) & 132.03 (2.35) & 0.92 (0.01) \\ 
2-step \verb+hw.glasso+ &  79.01 (0.08) & 98.08 (0.03) & 96.58 (0.22) & 100 (0) / 100 (0) & 150.42 (1.35) & 0.88 (0.01) \\ 
(known hubs) & \\
\hline
 &         & \multicolumn{3}{c}{$n=100, p=200$} \\
\verb+glasso+ &  64.25 (0.18) & 99.64 (0.01) & 63.59 (0.75) & 52.33 (1.66) / 100 (0) & 128.12 (3.53) & 1.50 (0.01) \\  
\verb+Ada-glasso+ & 70.13 (0.29) & 99.64 (0.03) & 80.73 (0.64) & 62.67 (2.69) / 100 (0) & 150.55 (6.42) & 1.08 (0.01) \\ 
\verb+SF+  & 67.81 (0.18) & 99.74 (0.01) & 81.60 (0.79) & 80.67 (1.85) / 100 (0) & 122.53 (2.33) & 1.15 (0.01) \\ 
\verb+hw.glasso+ & 69.47 (0.14) & 99.75 (0.01) & 86.45 (0.45) & 87.33 (1.63) / 100 (0) & 127.06 (1.74) & 0.96 (0.01) \\ 
2-step \verb+hw.glasso+ & 68.85 (0.27) & 99.41 (0.01) & 87.81 (1.27) & 88.33 (1.60) / 100 (0) & 191.34 (2.63) & 1.07 (0.01) \\ 
2-step \verb+hw.glasso+   & 70.79 (0.05) & 99.39 (0.01) & 96.86 (0.20) & 100 (0) / 100 (0) & 202.21 (2.38) & 1.05 (0.01) \\ 
(known hubs) & \\
\hline
 &         & \multicolumn{3}{c}{$n=250, p=500$} \\
\verb+glasso+ & 63.36(0.02) & 99.81(0.004) & 48.44(0.08) & 33.33(0)/100(0) & 367.27(5.06) & 1.88(0.01) \\ 
\verb+Ada-glasso+ & 64.68(0.04) & 99.97(0.00) & 52.43(0.14) & 33.33(0)/100(0) & 183.66(0.79) & 1.26(0.004) \\ 
\verb+SF+  & 63.53(0.02) & 99.96(0.001) & 49.08(0.06) & 33.33(0)/100(0) & 185.46(1.66) & 1.35(0.004) \\ 
\verb+hw.glasso+ & 63.50(0.02) & 99.97(0.001) & 49.05(0.09) & 33.33(0)/100(0) & 171.22(1.13) & 1.19(0.004) \\ 
2-step \verb+hw.glasso+ & 63.38(0.003) & 99.91(0.002) & 48.58(0.01) & 33.33(0)/100(0) & 240.18(2.47) & 1.26(0.004) \\ 
2-step \verb+hw.glasso+   & 72.21(0.10) & 99.89(0.002) & 80.81(0.38) & 99.00(0.57)/100(0) & 361.80(3.51) & 1.22(0.01) \\ 
(known hubs) & \\
\hline
 &         & \multicolumn{3}{c}{$n=1000, p=500$} \\
\verb+glasso+ &  65.28(0.09) & 99.63(0.01) & 54.85(0.28) & 33.33(0) /100(0) & 613.97(8.65) & 1.27(0.01) \\ 
\verb+Ada-glasso+ &  75.19(0.07) & 99.85(0.002) & 75.91(0.24) & 37.33(1.19)/100(0) & 437.80(2.24) & 0.74(0.002) \\ 
\verb+SF+  & 72.12(0.31) & 99.91(0.003) & 77.22(1.01) & 69.33(2.40)/100(0) & 330.97(6.51) & 0.80(0.004) \\ 
\verb+hw.glasso+ &  74.29(0.15) & 99.80(0.002) & 80.65(0.39) & 63.00(1.89)/100(0) & 492.86(4.11) & 0.67(0.003) \\ 
2-step \verb+hw.glasso+ & 72.04(0.36) & 99.87(0.002) & 77.78(1.23) & 63.67(1.90)/100(0) & 385.52(5.35) & 0.73(0.01) \\ 
2-step \verb+hw.glasso+   & 77.35(0.04) & 99.84(0.002) & 99.09(0.06) & 100(0)/100(0) & 476.43(2.35) & 0.64(0.002) \\ 
(known hubs) & \\
 \hline
\\
\end{tabular}
\caption{Means (and standard errors) of different performance measures over 100 
replications for the graphical lasso (glasso), adaptive graphical lasso (Ada-glasso), 
scale-free network approach (SF), hubs weighted graphical lasso (hw.glasso), 
and 2-step hw.glasso. } 
\label{tab:Simulation4}
\end{table}

\newpage

\end{document}